\documentclass[a4paper,12pt,leqno]{article}
\usepackage{geometry}
\usepackage{tikz-cd}\usepackage{bussproofs}
\usepackage{amssymb,amsmath}
\usepackage{proof}
\usepackage{psvectorian}
\usepackage{CJKutf8}
\usepackage[utf8]{inputenc}
\usepackage{graphicx}
\usepackage{graphicx}
\graphicspath{ {images/} }
\usepackage{qtree}
\usepackage{mathtools}
\usepackage{hyperref}
\usepackage{listings}
\usepackage{indentfirst}
\usepackage{titlesec}
\numberwithin{equation}{section}

\newtheorem{defn}[equation]{Definition}

\newtheorem{rem}[equation]{Remark}
\newtheorem{exm}[equation]{Example}

\newtheorem{notat}[equation]{Notation}
\newtheorem{newpar}[equation]{}

\newtheorem{xdefn}{Definition.}
\newtheorem{xproposition}{Proposition.}
\newtheorem{xcorollary}{Corollary.}
\newtheorem{xrem}{Remark.}
\newtheorem{xexm}{Example.}
\newtheorem{xlemma}{Lemma.}
\newtheorem{xtheorem}{Theorem.}
\newtheorem{xnotat}{Notation.}
\newtheorem{xnewpar}{\it}
\newtheorem{xproof}{{\it Proof. }}
\newtheorem{xproofof}{{\it Proof}}

\newenvironment{newparagraph*}[1]{\begin{xnewpar}\hspace*{-1.5mm}{#1}. \rm}{\end{xnewpar}}

\newenvironment{definition*}{\begin{xdefn}\em}{\end{xdefn}}
\newenvironment{remark*}{\begin{xrem}\em}{\end{xrem}}
\newenvironment{example*}{\begin{xexm}\em}{\end{xexm}}
\newenvironment{notation*}{\begin{xnotat}\em}{\end{xnotat}}
\newenvironment{proposition*}{\begin{xproposition}}{\end{xproposition}}
\newenvironment{corollary*}{\begin{xcorollary}}{\end{xcorollary}}
\newenvironment{lemma*}{\begin{xlemma}}{\end{xlemma}}
\newenvironment{theorem*}{\begin{xtheorem}}{\end{xtheorem}}

\titleformat*{\section}{\large\bfseries}

\begin{document}

\title{Hegel and Modern Topology}
\author{Clarence Protin}

\maketitle

\begin{abstract}
In this paper we sketch how some fundamental concepts of modern topology (as well as logic and category theory)  can be understood philosophically in the light of Hegel's Science Logic as well how modern topological concepts can provide concrete illustrations of many of the concepts and deductions that Hegel used.  Also these modern concepts can in turn be very powerful hermeneutic tools permitting a more rigorous and thorough grasp of Hegelian concepts. This paper can be seen as a continuation of our paper \cite{pro} where we argued that the prototypes of many fundamental notions of modern topology were already found in Aristotle's Physics.  More generally it is hoped that this note makes a case for the possibility of  a rigorous  enriching interaction and mutual support between philosophy on one hand and modern logic and mathematics on the other. This paper is obviously meant only as a preliminary sketch and to offer some motivation for exploring in a more detailed and thorough way the subjects discussed. \end{abstract}

\section{Introduction}

Why investigate such a topic as Hegel and Modern Topology? First of all, the author of this paper, acquainted with the major works of Hegel  and not being in any sense a 'Hegelian',  proposes a certain framework (the discussion of which must wait for a future paper) for assessing the significance of the Hegel in the history of philosophy with special emphasis on Hegel's relationship to ancient Hellenic, Hellenistic and Indian philosophy.  Doubtlessly such a framework provides some justification or motivation (and this is clear in the context of an Aristotelian\cite{pro} or Neoplatonist reading of Hegel)  for approaching the topic of Hegel's relationship to mathematics and more particularly the significance of Hegel for an interpretation of modern mathematics as well as in turn the significance of modern mathematics for an exegesis of Hegel. This is what Van Lambalgen\cite{pino} - in a Kantian context - aptly called the 'hermeneutical virtuous circle'.

Studying the mutual interpretation between Hegel's Science of Logic and modern mathematics is not a new proposal.  Hegelian concepts already appear in Lautmann's philosophy of mathematics.  The most advanced developments of modern geometry and topology are inseparable from Category Theory. This last provides a kind of abstract unifying conceptual architecture. It is plauisble that MacLane was conscious of the Hegelian undertones of his use of the tern 'subsumes' in his treatment of Kan Extensions in his book Categories for the Working Mathematicians. But the first explicit treatment of  the relationship of Hegel's Science of Logic and modern topology and category theory was given by William Lawvere (in the context of his philosophical approach to continuum mechanics and synthetic differential geometry) and there is major contemporary work (by the team responsible for the nLab website) in this direction which develops Lawvere's ideas.  This  work is focused on homotopy type theory, homotopical algebra and certain attempts at extending topos theory and synthetic differential geometry to a higher category theoretic framework, with notable applications to theoretical physics being carried out by Urs Schreiber.

Hegel's Logic (both in the Greater Science of Logic and the Encyclopedia Logic) does contain (principally in the Logic of Quantity and Measure)  an explicit treatment and critique of the foundations of arithmetic, geometry, infinitesimal analysis as well as applied mathematics. These portitions (or rather Hegel's extensive comments) have as a rule been dismissed as either inadequate or of merely historical interest (pertaining to a time before the establishment of the modern rigorous foundations for the Calculus). However around the beginning of the last decade interesting research has been carried out by Mascarenhas\cite{Fab} and Dimitrov\cite{dim} offering a more sympathetic reappraisal of Hegel's approach to mathematics in the Science of Logic and even a partial revindication of it in the light of subsequent developments in the foundations of mathematics.

The goal of this paper is to present a series of  sketches regarding how the virtuous hermeneutical circle could be carried out between Hegel's Science of Logic on one hand and modern topology on the other (as well as category theory and some aspects of categorical logic, type theory, mathematical physics  and theoretical computer science). 

In this paper we will only refer to basic 20th-century logic and mathematics such as taught in standard undergraduate degrees and in particular remain at the level of ordinary (1-)category theory except for some passages where we attempt to give an intuitive explanation for the basic ideas behind homotopy type theory and higher groupoids. 
But in the appendices we give general (philosophically motivated) introduction to the higher category theoretic and homotopy theoretic approaches found in nLab as well as discussing some more advanced aspects of topos theory and sheaf theory interesting in relationship to Hegel and also in their own right.

The presentation of the material here is obviously very sketchy and its primary aim is to inspire and motivate future research in this area.

\section{General considerations on Hegel's Logic}

Hegel expresses apparently negative views of mathematics in his published works (specially in the Science of Logic). However we need to balance our judgment with the following extract from his letters (see Wallace's translation of the Encyclopedia Logic, xiv-xv) which is very revealing:\\

\emph{'I am a schoolmaster who has to teach philosophy, who, possibly for that reason, believes that philosophy like geometry is teachable, and must no less than geometry have a regular structure. But again, a knowledge of the facts in geometry and philosophy is one thing, and the mathematical or philosophical talent which procreates and discovers is another : my province is to discover that scientific form, or to aid in the formation of it'.}

\emph{' You know that I have had too much to do not merely with ancient literature, but even with mathematics, latterly with the higher analysis, differential calculus, chemistry, to let myself be taken in by the humbug of Naturphilosophie, philosophising without knowledge of fact and by mere force of imagination, and treating mere fancies, even imbecile fancies, as Ideas.'}\\

Our view is that Hegel is essentially concerned with the phenomenology of consciousness (including forms of social, historical, cultural and religious consciousness) and that the Science of Logic represent an attempt at abstracting the structures and processes found therein. True, there is a sense in which consciousness is left behind at the beginning of the Logic, that it proceeds from where the Phenomenology of Spirit left off, that even pure being is incomprehensible by consciousness.  Be this as it may, it is still the case, so we contend, that consciousness is the ultimate source of the Logic. Consider this passage from the Encyclopedia Logic:\\

\emph{In other words, every man, when he thinks and considers his thoughts, will discover by the experience of his consciousness that they possess the character of universality as well as the other aspects of thought to be afterwards enumerated. We assume of course that his powers of attention and abstraction have undergone a previous training, enabling him to observe correctly the evidence of his consciousness and his conceptions.}\\

Also  we read in the Lessons in the History of Philosophy:\\

\emph{The two formal moments in this sceptical culture are firstly the power of consciousness to go back from itself, and to take as its object the whole that is present, itself and its operation included. The second moment is to grasp the form in which a proposition, with whose content our consciousness is in any way occupied, exists. An undeveloped consciousness, on the other hand, usually knows nothing of what is present in addition to the content}.\\

The following appear to us to be some of the fundamental processes of consciousness which also correspond to certain fundamental mathematical structures or constructions:

\begin{itemize}
\item Synthesis. This is reflected in the gluing of geometric objects, various kinds of covers and specially the sheaf condition. 

\item Self-reflection. A system which can represent (partially at least) higher order aspects of itself within itself.  A prime example is the original synthetic unity of apperception:  \emph{I know that I am knowing}. This is reflected in recursive and inductive definitions. 

\item Return-to-self. This is related to the negation of the negation (double negation). An example is found in temporality: something must pass away to reveal its essence - to ti ên einai, literally 'the what it was to be'. 

\end{itemize}

Let us recall the basic trinary structure of Hegel's Logic (see \cite{car} for a deeper discussion). The Understanding seizes a concept as an immediate whole, but then Dialectical reason remembers the complex process of its genesis and its inward relations and tensions (that it, it becomes a category of objects and morphisms).  Speculative reason then replies to Dialectical reason that these relations actually flow in both ways rendering the category again a coherent whole - maybe we can view this as realizing the quotient construction in which we find formal inverses to the relations brought forward by Dialectical reason. Or it is like taking an integral of a certain field density over all possible paths between two points. Thus Speculative reason subsumes into bilateral  isomorphisms these directed unilateral relations and again brings peace to the concept, then to be seized again by the Understanding at a higher level.

In sections 2-9 we only consider elementary mathematics  while in the remaining sections we consider advanced category theory.

\section{Hegelian Understanding as necessary amnesia and unanalyzable identification: univalence and coherence}

In an upcoming paper by Carlos Lobo on Husserl's project of the reform of logic it is spoken of \emph{la strate-substrat dont elle (la logique formelle) suppose tacitment l'inanalysabilité comme condition de son propre fonctionnement}, which evidently recall's Girard's theory of the "blind-spot".

The idea of formal logic necessarily supposing a non-analyzable substrate for its own functioning is certainly a thought-provoking proposal which would point to, in Kantian terms,  a necessary synthetic a priori principle or condition of human cognition.

But here we will only briefly discuss an interpretation of a such a principle in terms of Hegelian dialectic as well as some striking illustrations (or confirmations) of this interpretations offered by type theory and category theory.

The problem of the individuality of formal objects cannot be separated from problems of equality, 
identity, extensionality vs. intensionality - nor indeed from the historically rooted ontological-logical 
problems of genus, species, difference and what could be the minimal elements (for the predication 
relation) in the hierarchy of species and genera. In Aristotle the minimal element in the category of 
substance is the individual substance.  In \cite{pro2} we inquire about the analogous minimal elements for the other categories such as quantity, quality or 
relation. We can ask both what constitutes the individuality of a formal object and what an individual 
could be in a class of formal objects. A fundamental constitutive act in question seems to be: 
identification under isomorphism (such as when we say let $\mathbb{N}$ be *the* set of natural numbers...).  
Identification under isomorphism is a manifestation of such an "unanalysability", a kind of positing of 
indivisibility, indecomposability, indiscernibility which calls to mind an indivisible continuum (cf. 
Aristotle's\emph{genos hôs hulê},  genus is like matter). 

In  Martin-Löf type theory - the extension of ML type theory to 
homotopy type theory. In HTT logical objects, types/propositions/sets, take on the nature of higher 
groupoids endowed intrinsically with a qualitative-deformative notion of equivalence.
With the work of Voevodsky and HTT a fundamental logico-transcendental supposition of mathematics 
was brought to light: the univalence axiom.  This axiom can be paraphrased as follows: what makes 
mathematics possible to be formalized (or more naturally formalizable) on a ML type-theoretic basis is the necessary assumption that 
we can treat equivalent (i.e. isomorphic) objects as being equal in the strict sense (i.e. their differences 
are unanalyzable in the system). This transcendental-logical axiom is expressed formally within HTT 
itself in the form: equivalence is equivalent to equality.

 For Hegel the triadic cyclic ascending process which passes through the understanding to dialectical reasoning to speculative reasoning back to the understanding at a higher level is seen as a necessary framework for knowledge. We could interpret the Hegelian understanding as involving precisely this non-analyzable substrate and identification (merging, synthesis) necessary for knowledge to progress. In fact, the understanding involves a voluntary amnesia regarding an object's process of coming to be,  an attempted static synthesis wherein all the past varied dynamical genetic aspects are forgotten and merged into a tranquil unanalyzable unity, a kind of ”premature synthesis” which according to Agnieszka Matylda Schlichtinger can be interpreted as a stage in the mergence of order out of chaos. The Hegelian understanding can be construed as an imposed condition for static objectivity.

The appearance of the understanding on the stage is analogous to transforming a system with memory into system without memory in which the state of the system only depends on the initial and final input and not on the path connecting one to the other. We can even imagine that the states of the memory-less system are obtained by collapsing or identifying all such paths between two inputs. In a monoidal category without the coherence conditions if we take any finite collection of objects $A_1,…,A_n$ then there are finitely many different objects we can form by $\otimes$- tensoring these objects to form new objects X. If we consider two objects X and Y formed in this way then in general there is more than one path, more than one composite isomorphism involving tensorings of id and associators, linking X to Y. The coherence condition states precisely that all these composite isomorphisms must be the same (it has been shown that the famous pentagon condition is sufficient to ensure this for arbitrary X and Y): this is clearly an expression of the lack of memory of the system (the isomorphism must depend only on X and Y) and of a merging, synthesis, unification and voluntary amnesia regarding the difference between possible composite isomorphisms linking the two objects. The coherence condition expresses a necessary condition for the category to be able to model $\otimes$-constructed objects in an consistent way via isomorphisms, we could almost say, that the coherence condition guarantees that there is a consistent internal logic of the category.

\section{First overview of the Logic of Being}

The first section of Hegel's Logic, the Logic of Being,  involves an abstract exploration of the concept of  space and the differential geometry used in mathematical physics. Hegel, who taught differential calculus both at high school and  university, dedicates a long section of the Logic to the infinitesimal calculus, the notions of which illuminate many other passages of the Logic. 

In the section on Being, in the parts involving the One and the Other, Limit, the Infinite, the One, and the successive transitions involved - all this suggests a connection to a universal construction in mathematics: that of a completion. We have a finite structure which is limited, incomplete in some determinate way. For instance a field may not be algebraically closed. Something is missing, something extrinsic that needs to become intrinsic (cf. the process of finite algebraic extensions, adjoining roots, etc.).  Yet the finite completion, the finite adjoining of the missing aspects generates in turn its own incompleteness at the next level, and so on. This situation is only overcome by assimilation, by incorporation of this extension process, this outflow, into the structure itself - now invariant under this finite 'passing to the limit'. This 'incorporation' is itself a kind of  higher 'limit', a limit of limits, so to speak (which has a universal property, or a minimality property). We find similar processes in measure theory (in the construction of the Borel $\sigma$-algebra, regular and complete measures) and analysis (completion and approximation theorems). Hegel states in the Encyclopedia Logic that the determinations of Being are exterior to each other, their process is the passing into another.

The transition for Being-for-self to Quantity recalls the proof of Urysohn's lemma. The main ingredient is the fact that in a locally compact Hausdorff space $X$ if we have a compact set $K$ contained in an open set $U$ then we can find an open set $V$ such that $K \subset V \subset \overline{V} \subset U$. This expresses $K$'s tendency to overcome its limitation relative to $U$. In the proof of Urysohn's lemma this process of overcoming the limit is multiplied to countably infinite sets indexed by the rationals and then completed into a continuous function, into continuous quantity. Somehow this relationship to the real numbers is implicit in the abstract concept of a locally compact Hausdorff space. There are analogous situations in which valuations, measures, uniparametric semigroups of automorphisms, emerge.  For example for continuous functions with compact support  Borel measures emerge from bounded positive linear functionals on spaces of such functions. 

Pure Quantity or indifferent Quantity (apeiron) is topological, corresponding to what is implicitly determined by a variable, a sequence, a net (cf. the discussion in \cite{pro}). Attraction is aggregation, cohesion, expressed by the open sets of a topology.  Repulsion is discreteness . The duality and mutual transition between the continuous and discrete is reflected in the fact that geometric structures give rise to algebraic structures and these in turn have their geometric realizations. 

In the initial section of the Greater Logic on Measure there is a suggestion that Measure is a kind of conservation law; that which remains invariant under actual, possible or conceivable change and variation. Measure is related to mathematical physics; Hegel states that a Measure  is not pure Quantity but Quantity in relationship to something exterior to itself. Hegel also raises the question of natural vs. artificial units of measure. He introduces gradualness (i.e. continuous change in quantity with surprising jumps in quality) and mentions the Sorites paradox. Hegel's interpretation might be couched in terms of locally constant sheaves (equivalent to covering spaces):  local homogeneity and uniformity is reconciled with global non-triviality. Hegel's knotted line foreshadowed the modern notion constructible sheaf\cite{Kashiwara}, the 'knots' corresponding to the stratification (cf. singular points of a variety, bifurcation set, etc.).

The section on Measure concerns also how systems combine into larger systems or split again, etc. Hegel uses chemistry and thermodynamics. It suggests the idea of program or function which can be integrated with different inputs producing different behaviors yet still maintaining its essence.

\section{More on the Logic of Being}

The something and the other. In the Science of Logic Hegel thinks of the idea of two things bearing the same relation to each other and indistinguishable in any other way. Now a good illustration of this situation in found in the two possible orientations of a vector space. Each bears the same relation to the other and there is absolutely no way to uniquely identify one of them in distinction to the other.  In the same place Hegel talks about the meaning of proper names and states that they have none !

Hegel's strategy. All of Hegel is based on the dynamics and structure of consciousness as it reveals itself to itself. But equally important it is all based on positing this structure of consciousness  to be essentially objective and not subjective.  Thus 'being a subject' is seen as a stage in the development of the object while mere  'subjectivity' is seen as a failure and partiality of objective consciousness in not living up to its full objectivity or actuality (or knowledge and realization thereof).

Thus ideality and infinity are the key structure of an object which has developed to a point of being a subject. 

Finitude and limit:  a closed set.  Yet the boundary also defines what the set is not (cf. discussion on Heyting algebra of open sets, etc.). Thus the boundary contains implicitly its own negation.

Limitation : an open set (also analytic continuation). Limit outside itself.

The ought... manifestation of the germ of a space through a given open set representation of the equivalence class. Any particular one is insufficient and can be replaced by another which is also insufficient.

False infinity: taking simply the set (or diagram) of such  open set representatives.  True infinity taking the limit (equivalence class) or limit (in the categorical sense)..

Degrees of interpenetration and transparency between the individual and the universal (and between self-consciousness and essence) in the Phenomenology of Spirit.  From the rudimentary form in the ethical life to: the universal is in the individual, the individual in the universal and the universal is in the relation between individuals and the relation between individuals is in the universal (mutual confession and forgiveness). 

Milne's notes on Class Field Theory has a quote of Chevalley (1940) which seems couched in Hegelian terms: the Abelian extensions of a number field are determined by the number field itself, the number field  'contains within itself the elements of its own transcendence'.

And yet : and yet it is really basic to use ontology and logic inspired by (or abstracted from) biology and psychology  - that is, a kind of general systems theory and universal logic based fundamentally on biology and consciousness. The fact that mathematics exhibits philosophical analogy with biology suggests that biological ontology and logic is fundamental. René Thom already did this by suggesting that concepts have analogies with living organisms.

'Pure' category theory is very poor in results, like Hegel's poor, empty, abstract universal.  But once we consider the particular, like Abelian or Monoidal categories, things become very rich.

Extensive magnitude: the category including in itself objects and arrows. Intensive magnitude, one object signaled out, hence only meaningful in the context of a category.  Infinite progression:  a diagram in a category. Ratio the limit. For infinitesimals the limit for the category of open sets does not exist but for the functor it does.  Ratio is also the quotient construction in which representatives are sublated moments.

\section{On Quantity}

The indifference of Quality expresses an infinitesimal local Quality around a point. The idea of an open set in which the actual size is completely irrelevant, the local homogeneous quality. This can be understood in terms  of a basis of a topology $\tau$: a family of open sets such that each open set $U$ containing a point $x$ contains an element of the basis (containing $x$). The ordinary topology of $\mathbb{R}^n$ has a basis given by open balls of radius $\frac{1}{n}$. The Hegelian 'indifference' is expressed by the fact that we can remove any finite (or even some infinite sequences) of such balls from the basis without affecting the topology defined.

The definition of topology is inseparable from the concept of continuity. And yet the ordinary topology has an element of discreteness, or rather, can be generated by a discrete sequence.

The concept of continuity, in our view, is  tied to that of sheaf, in particular to the concepts of gluing or synthesis. In its most rudimentary form it is the condition that any open set must be able to be written as a union $\bigcup_{i\in I} B_i$ of open sets of a basis. 

Consider a category $C$ of open sets. And let $D: I \rightarrow C$ be a diagram, for instance, consisting of open balls containing a given point.  In this case the point is a 'limit' of the diagram in some sense, but not in the category $C$. But we can define the notion of $c$-limit - the idea here is that the limit point is defined implicitly through the complement of this maximal open set.

Given a diagram $D: I \rightarrow C$ for $I$ a directed set, we say that an open set $V$ in $C$ is a $c$-limit of $D$ if i) for any object $i$ in $I$ we have that $F(j) \cap V \neq 0$ for some $j \prec i$ and any open set $V'$ satisfying this condition is contained in $V$. For example, for a sequence of balls centered at the origin the complement of the origin is the $c$-limit.

Hegel's treatment of repulsion and attraction  recalls the separation axioms in general topology. And the concept of connectivity is also directly involved with separation.  The extremes are that of the coarse and discrete topologies. In the coarse topology everything is merged together, we cannot separate things. In the discrete topology not only can we separate everything but each point is its own separation (this atomization is a form of alienation that features in many crucial moments of the Phenomenology of Spirit).

Hegel's 'Units' are like the proximate species of a genus. In Aristotle genus is a continuum, like matter. The species carve up the genus into pieces (which themselves become further genera). But all species form a spectrum, a continuum in sharing the quality of the proximate genus.  All open sets (at least basic ones) should share the same quality. The topology is the quality. In the sheafication construction of a presheaf, a section of the constructed sheaf is essentially a unity, a gluing of local sections, a kind of averaging out into a locally homogeneous quality.

Hegel's treatment of number in the section of Quantum is noteworthy. Given a non-empty finite set $S$ for every $s\in S$ we have $S = S' \cup \{s\}$ with $S' \subset S$. This is what Hegel seems to call $s$ being a limit of $S$. It matters not which $s$ is used. Hegel's whole approach suggests an intuitive continuum-based phenomenological view of number, numbers as connected objects which expresses the internal homogeneity and indifference between objects.

Different total orderings of a finite set are equivalent.  For each we have  $S = \{n \in N: n \leq max\,S\}$. Orderings are external, imposed.  Recall also the external relation of Frege's equinumerous relation. A number is the set of all sets equinumerous to a given set (the extension of the number)  and again the precise choice of this set is irrelevant.

Hegel was certainly prescient about the relationship between geometry and arithmetic, the discrete and the continuous. Modern mathematics has confirmed this.  Numbers do not add themselves, they are added from outside. Arithmetical operations are imposed from without, much like the category theoretic view of algebraic structures (for instance a natural number object).  We could also wonder if intensive quantity corresponds to viewing the real numbers as a set with an order. Hegel seems to have anticipated in this discussion the axiom of induction in Peano Arithmetic.

Homotopy (in algebraic topology) is the passage of Quantity into Quality. It is a changing of shape and size which preserves - and thus defines - a certain quality.

A simplicial set is an abstraction of a topological space, it is a category theoretic abstraction of a geometric form (i.e. a polytope). But it is a geometric form which contains within itself the process of its own genesis or assemblage (cf. the geometric realisation functor associating a topological space to each simplicial set). 

The category of sets $Set$ is an exterior, abstract, discrete concept (Hegel's concept of number seems to be very set theoretic), but the category of simplicial sets $sSet$ represents a greater cohesion between parts and qualitative structure and determination. The morphisms between objects in a simplicially enriched category instead of being a mere set become a space.

\section{The transition from Being through Quality to Quantity}

Consider a category $C$ with only one object $A$.   Is $A$ then an initial or terminal object ?  Not necessarily either because $hom(A,A)$ can contain other morphisms besides $id_A$ (i.e. we might say it has internal symmetries, it is not absolutely simple).   If $A$ is a terminal object then $hom(A,A) = \{id_A\}$,  but then $A$ is also an initial object.   And if $A$ is an initial object then $hom(A,A)$ again must be $\{id_A\}$ and so $A$ is also a terminal object.  This is one  (very simple) illustration  of the passage between Being and Nothing and Nothing and Being.

If $C$ is a category with two objects $A$ and $B$ and if $A$ is initial and $B$ is terminal (there is a canonical $t: A \rightarrow B$) then if there is a morphism $f: B \rightarrow A$ then $A$ and $B$ must be isomorphic. This is a simple prototype of the quotient construction.  But let $C$ have three objects, initial $0$, terminal $1$ and $A$.   Then it follows from these hypotheses that the canonical  $0  \rightarrow  1$     must factor a unique way as   $0 \rightarrow A \rightarrow 1$.  Now $A$ is allowed to have internal symmetries ($hom(A,A)$ need not be $\{id_A\}$).   And postulating a $1 \rightarrow A$ means that "$A$ has an element".  Thus $A$ is now a Something.

On a very simple level in a category we find already the embodiments of the basic Hegelian concepts: Being-in-self,  Being-for-other and Being-for-self.    Objects, considered in themselves, without arrows, without relations, are very poor and empty.  The only morphism an object $A$ must have is $id_A$. We could almost say that $id_A$ represents $A$ (in some formulations of a category as a partial monoid we can even identify objects with these morphisms).  $id_A$ represents the bare minimal immediacy of a quality (or being at this stage).  This is being-in-self.   Now the existence of a morphism $f: B \rightarrow A$ represents a relation between $A$ and $B$, it represents $A$'s being-for-other,  $A$'s being relative to $B$.  Likewise $g : A \rightarrow B$ represents $B$'s being relative to $A$. The composition of these two   $f \circ g : A \rightarrow A$   is in general non-trivial (distinct  from $id_A$) abd represents being-for-self, a more complex self-relation which is mediated by relations to the Other.

Much topological intuition and anticipation of modern geometry and logic seems to be found in Hegel's Science of Logic.
For instance Hegel writes: the boundary of a thing both belongs and does not belong to it (we discuss the relationship to the Heyting algebra of open topological sets \cite{pro}  and the logical significance of the boundary in a topological sense).   Points are seen both as limits of lines and generators (in the algebraic and indeed Lie group sense). The pure something in Determinate Being (when we abstract it from its Limit) is likened to a local homogenous quality (open set). The limit is seen as a interface with the environment.  Constitution (for instance rational thought is called the constitution of man) implies a barrier to be overcome and already a certain rudimentary invariance under change. Hegel even says points, lines, surfaces have a "dialectical tension".

Infinity is the negation of the negation (the connection to the double negation topology in topos theory\cite{Moerdijk} is worth exploring - it is connected to density and indeed probably to self-similarity and scaling). But there is a simple illustration that might be of some interest.  Consider the Heyting algebra of open sets of a topological space (which is in general not a Boolean Algebra so we do not have $A \vee \neg A$ or equivalently  $\neg \neg A \rightarrow A)$. The negation of an open set  $A$,  $\neg A$, is defined as the interior of its complement  ($int\, A^c$).   We have $A \subset \neg\neg A$.  But not in  general an equality.    $\neg\neg A$ here represents the Hegelian infinite, which is essential dynamic, it is about evolving, overcoming barriers.   Now look at the open set $(0,1) \cup (1,2)$.   If we think of this open interval as having a living, flowing nature, what would the natural tendency be of this open set ? It would be to fill in the whole (the gap $\{1\}$), to cohere and merge into $(0,2)$.   But if we calculate  $\neg \neg ((0,1) \cup (1,2))$  we get precisely $(0,2)$.   The open sets which are invariant for the operation of double negation  are the regular open sets.

There is more that could be said about the One in terms of infinitesimals, stalks of sheaves and germs of functions and their unfoldings. In particular about how the One expels and generates other Ones. Hegel indicates that he is thinking of Leibniz's monads, of  infinitesimal points.

On an elementary level we can think of the One as a fundamental system of neighbourhoods around a point or equivalently as the germ of a point $p$ which is an equivalence class in a quotient construction. Thus it must have representatives actual open sets $U$ containing $p$. But then $U$ itself generates (or represents) germs for all the points $q$ in $U$.  Thus each germ will generate a multiplicity of others.    A slightly more advanced way of looking at this is in terms of germs of functions (around a singularity)  Each germ of a function determines a universal unfolding which gives rised to qualitative local discontinuities around the point (a catastrophe set).  Thus the homogeneous strata = attractions and the germs of submanifolds representing the catastrophe sets represent repulsion. In catastrophe theory one speaks of a "conflict" of attractors.

The idea of a germ of a point $p$ is quite a  good illustration of Being-for-self (the One).  The germ of a point is defined as an equivalence class (hence via a quotient construction) of all open neighbourhoods $U$ containing $p$. Thus $U$ and $V$ determine the same germ if $p$ is in $U$ and $V$ and $V$ is contained in $U$.   For $p$ to be a germ it needs to enter into relationship with the ambient space, it needs the topology of the space. It cannot be a germ merely in itself.   Thus to define the germ we need some open set U (the germ is for-another).  This $U$ is posited. But at the same time we can deny and discard this $U$ and replace it by a different one $V$ as long as it is equivalent to $U$ (for instance $V$ is contained in $U$ and contains $p$). This is the moment of negation. That is, the germ is not this, not that, etc. It has abrogated  complete determination via any specific open set.  The whole net or local basis around $p$ represents the 'bad infinity'.   The quotient construction that defines the germs represents finally the return-to-self and a unity.

This germ seen as a One has a "dialectical tension" which we can associate topologically with the idea of instability.  It has to unfold and stabilize in the ambient space but by doing so it reveals its internal contradictions and movement of attraction and repulsion. Maybe there is a connection to symmetry breaking in theoretical physics.

Besides the quotient construction another fundamental construction is that of completion (but completions and limits in general are often defined using quotient constructions).  Take the category of open sets of a topological space. Then we can define a category theoretic diagram  $U_1 \rightarrow U_2\rightarrow U_3 \rightarrow ....\rightarrow...$ (equivalent to a net, or local basis) which represents a germ at a given point. But, as we already remarked, the limit in the category-theoretic sense does not in general exist. Thus we need to complete the topology by adding infinitesimal neighbourhoods. So if $G(p)$ represents the germ-object and $U$ some neighbourhood we now have canonical morphisms  $G(p) \rightarrow U \rightarrow G(p)$ which represents return-to-self.

These germs or infinitesimals have both the nature of discreteness (they are point-like in the sense that for a "decent" topology two distinct infinitesimals are represented by disjoint open neighbourhoods around their respective points) and continuity. Open sets are obtained from infinitesimals by "integration", a kind of synthesis or gluing together of infinitesimals (cf. the Leibniz notation for the integral).  Thus the category of Quantity is born.  This is also related to the sheafication construction discussed above.

\section{Double-negation, sheaves and the Logic of Being}

As we shall see, double-negation plays an important role in the logic of Essence. However as a continuation of the previous section we offer some considerations on double-negation in topos theory and how it evidently reflects  many aspects of the Logic of Being.

We believe that Lawvere's original theory of elementary toposes and its subsequent development is important and interesting. The most important object in an elementary topos is the subobject classifier $\Omega$. If we think of the topos as a cell then $\Omega$ is a kind of nucleus. An elementary topos is simply a model of mathematics which as a foundation is vastly superior and more cogent than ZFC. To us the most important concept in topos theory is that of the Lawvere-Tierney "topology", which is just a morphism $j : \Omega \rightarrow \Omega$ satisfying three simple "modal" or "closure-operator" type axioms. A very important instance of $j$ is given by double-negation $\neg\neg : \Omega \rightarrow \Omega$, in which we consider the (internal) Heyting algebra structure on $\Omega$. The morphism $j$ then determines a localization of the topos, a new subtopos of the original topos called the topos of $j$-sheaves. 

A central philosophical problem of topos theory is understanding the meaning of the Lawvere-Tierney "topology" and its associated topos of $j$-sheaves as well as the special role of the $\neg\neg$-topology (why is it not abuse to call $j$ in this case a "topology"?).  A key to this is to see how the theory above abstracts the concrete case of presheaves and sheaves. A "sieve" is a curious concept. Think of a set $S$ of open sets (conceived as "cover") in some space $X$. Then take the minimal extension of $S'$ of $S$ under the condition that $U \in S$ and $V\subset U $ implies $V \in S $. Then we have a sieve (generated by $S$). We could rephrase the condition as $W \cap U \in S$ for any $U \in S$ and open set $W$. That is $S'$ is the $\wedge$-ideal generated by $S$. The we have the obvious notion of a principal idea generated by $O$ (called a principal sieve).

 For the presheaf topos on a topological space $X$ we have that the presheaf $\Omega$ associates to each open set $U$ the set of all sieves on $U$. So $\Omega$ is a kind parametric version of local truth values. On the presheaf topos an important example of $j$ is the functor that associates to each sieve $S$ (on a $U$) the principal sieve determined by "what $S$ covers". In the words of Moerdijk and Maclane "What counts is what gets covered". Thus in this case $j$ is nothing more than a kind of parametric generalized union. Logically it is expressing "if something is locally true then it is globally true". That is the subjobject of $\Omega$ determined by $j$ consists in those sieves which are invariant under generalized union, situations in which if something is locally true then it is globally true. This fails for instance for the presheaf of constant functions. This $j$ (which we should call the union topology) seems to be intuitively clear but we still need to understand better why in the topos of continuous functions on  topological space $X$ the internal logic of the topos proves that all functions are continuous.

Thus for the topos of presheaves on a topological space the subobject classifier on an open set $U$ yields all the sieves on $U$ while for the topos of sheaves it yields all the open sets of $U$ (or equivalently the set of principal sieves on $U$. 

But a central problem of topos theory is understanding other $j$s such as the double-negation topology (and its "$j$-sheaves"), in particular as an abstraction of the topological sheaf case. We have written something about this in our "Hegel and Modern Topology". The double negation topology is all about "density" while the union topology is about local-global coherence. Could we associate the union topology with exponentials in linear logic ?

What does the double-negation topology look like for the presheaf topos on a category $C$ ? And for the topos of sheaves over a topological space ? In this case it appears to be simple: an open set $V$ in $\Omega(U)$ is taken by $j = \neg\neg$ to $int(\overline{V})$. Thus the topology subobject $J$ is just the set of regular open sets in $\Omega(U)$. What are the $\neg\neg$-sheaves in this case ? The idea is like the passage from holomorphic to meromorphic functions (or the identification of measurable functions different on measure zero sets). The resulting sheaves can be considered as defined up to empty interior closed subsets. In a way, the $\neg\neg$-topology introduces closed boundaries, qualitative differences between different regions.

The fundamental definition of $j$-sheaf, an object satisfying $hom(E,A) \cong hom(B,A)$ for a $j$-dense subobject $E$ of a an object $B$ and considering the restriction morphism, can be understood as follows. The associated closure operator is a kind of transcendence, expansion beyond itself, saturation, perfection. The condition then picks objects for which their effect on a subobject already contains their effect on that subobject's transcendent expansion (closure). This means that the object satisfying the condition already contains within itself its own self-transcendence.

Put in another way, sheaves for the dense topology are ideal cohesive structures which on are object (open set) are determined by coverings whose union is not equal but only dense in that object. These sheaves allow a form of completion, limit, coherent transcendence. They allow one to construct an object based on coherent data which is yet still incomplete.

\section{Hegel's Doctrine of Measure and Thom's Catastrophe Theory}

The section on Measure in the Science of Logic is highly suggestive from the point of view of modern topology. We argue that modern differential topology and singularity theory represents a development of  key intuitions and concepts found therein.

The relationship between Quality and Quantity in Measure recalls  not only that of a physical field spread out through space-time, but the general concept of a sheaf. The base category is often space, or space-time, or some control space, while the target category is phenomenological, qualitative although having in turn a quantitative aspect.  The base is extensive, the fiber intensive. The further development in this section of the Logic - and most specially the concept of Nodal Line - appears to us as a prototype of Thom's Catastrophe Theory expounded, for instance, in his book Mathematical Models of Morphogenesis.  Hegel offers an important insight on the essential co-dependence between the qualitative and quantitative which has consequences for applied mathematics.  Strictly determined systems have to be studied qualitatively. But qualitative studies are only meaningful and fruitful in the context of a specific determined object. Thus we must employ qualitative concepts to talk about a sculpture (for instance locally the quality of the curvature, the topology of the surface levels) but a merely  global qualitative (i.e. topological) study would identify it with a ball. In order words, the topological exists to study the metrical and the metrical has to be studied using the topological. The same goes for differential equations and dynamical systems. As Hegel wrote: the Quantum is Quantity (in the sense of indifferent quantity) endowed with Quality.

Hegel's discussion on the Sorites  besides being related, as we mentioned, to locally constant sheaves, also invokes phase transitions in modern physics.  And is not general relativity a partial realization of Hegel's project of deriving physics entirely from the properties of space and time ?  Hegel mentions the influence on the observer on the experiment (in his discussion of thermodynamics) .  The  transition to Essence suggests that while applied mathematics is legitimate in its own domain and perspective, deeper scientific knowledge has to be found elsewhere - but this latter (whatever it be) must in all cases involve passing through a thorough understanding of the former.

Hegel's concept of Type suggests the image of a continuous map $F: U \rightarrow \Phi$ for some phenomenological space $\Phi$ where $U$ represents the allowed variations of the independent parameter. Besides the connection to Thom there is also the more general connection to moduli spaces.

\section{On Quantity}

The indifference of quality which expresses an infinitesimal local quality around a point. The idea of an open set in which the actual size is completely irrelevant. The local homogenous quality. This can be understood in terms of one form of the definition of basis of a topology $\tau$.  It is a family of open sets such that each open set $U$ containing a point $x$ contains an element of the basis containing that point. The ordinary topology of $\mathbb{R}^n$ has a basis given by open balls of radius $\frac{1}{n}$. The Hegelian 'indifference' is expressed by the fact that we can remove any finite (or even some infinite sequences) of such balls form the basis without affecting the topology defined by the basis.

The definition of topology is inseparable from the concept of continuity. And yet from the exposition above we see that even the concept of the ordinary topology has an element of discreteness, or rather, can be generated by a discrete sequence.

The concept of continuity, in our opinion, is intimately tied to that of sheaf, in particular with that of gluing or synthesis. In its most rudimentary form it is the condition that any open set must be able to be written as a union $\bigcup_{i\in I} B_i$ of open sets in a basis. 

Consider a category $C$ of open sets. And let $D: I \rightarrow C$ be a diagram, for instance, balls containing a given point.  In this case the point is a 'limit' of the diagram in some sense, but not in the category $C$. However we can define the notion of $c$-limit.

Given a diagram $D: I \rightarrow C$ for $I$ a directed set, we say that an open set $V$ in $C$ is a $c$-limit of $D$ if i) for all all object $i$ in $I$ we have that $F(i) \cap V \neq 0$ and any open set $V'$ satisfying this condition is contained in $V$. For example, for a sequence of balls centered at the origin the complement of the origin is the $c$-limit.

Separation and attraction.  Recall the separation axioms in general topology. And the concept of connectivity is also directly involved with separation.  The extremes of the coarse and discrete topology. In the coarse topology everything is merged together, we cannot separate things. In the discrete topology not only can we separate everything but each point is its own separation.

Hegel's 'units' are like the proximate species of a genus. In Aristotle genus is a continuum, like matter. The species carve up the genus into pieces (which themselves become further genera). But all species form a spectrum, a continuum in sharing the quality of the proximate genus.  We can think of a connected groupoid of connected groupoids. Hegel's unit is like a object in a category. The category of open sets we discussed above. All open sets (at least basic ones) should share the same quality. The category, the topology, is the quality.

Thus we should view quantity as a Kantian process of subjective synthesis - which we already interpreted in terms of the sheaf-gluing process. In the sheafication process we see that a section is essential a unity, a gluing of local section.

Hegel's treatment of number in the section of Quantum is curious. Indeed given a non-empty finite set $S$ for every $s\in S$ we have $S = S' \cup \{s\}$ with $S' \subset S$. This is what Hegel seems to call $s$ being a limit of $S$. It matters not which $s$ is used. Hegel's whole approach suggests an intuitive continuum-based phenomenological view of number, numbers as connected groupoids which expresses the internal homogeneity and indifference between objects.

Different total orderings of a finite set equivalent.  For each $S = \{n \in N: n \leq max S\}$. There are external, imposed.  Also external relation of Frege's equinumerous relation. A number is the set of all sets equinumerous to a given set (the extension of the number ?) ... and again this set is irrelevant. An equivalence class. A groupoid.

Hegel is certainly prescient about the relationship between geometry and arithmetic, the discrete and the continuous. Modern mathematics has confirmed this.  Numbers do not add themselves, they are added from outside. Arithmetical operations are imposed from without. This is just the category theoretic view of algebraic objects in a category or a natural number object ! 

Hegel seems to have anticipated induction. Intensive quantity, degree: here we have number as the real numbers, as a set with an order. Now limit is exterior. 

\section{More on Quantity and Measure}

The  connections between Hegel and Kant, Fichte, Schelling, Leibniz and Spinoza have been substantially explored.  However there are lesser known and perhaps even more philosophically interesting connections to Sextus Empiricus (as evidenced
in Hegel's Lectures on the History of Philosophy and of course in the Phenomenology itself), Hume and Proclus. The same goes for the various currents of ancient Indian philosophy (themselves with interesting connections to Pyrrhonism and Hume).  A topic to be explored is Hegel's relationship to Proclus' theory of dialectics and philosophy of mathematics in the Commentary of the the first book of the Elements.
 
  A major work relevant for studying Quantity is the  PhD thesis of F.M.  Nolasco\cite{Fab}.  This thesis suggests the central role of the question of the \emph{algebraization of the calculus} in Hegel's critique of the mathematics of his time. 

Concerning the section on quantity and the theory of infinitesimals in particular,  a central role is played  by what might be considered  the algebraic expression of infinitesimals (and continuous quantity): attempts to capture infinitesimals via sequences, power series and in particular convergent power series. And the same goes for the related concepts of point, local neighborhood, etc. 
Perhaps the most interesting embodiment of these concepts is found in algebraic geometry. The idea of completion of a local ring relative to its maximal ideal, which yields a finer or smaller kind of locality (expressed  for algebraic varieties by a quotient of a power series ring) then just  the local ring.  In the standard treatment of completion (relative to a topology which comes from a filtration defined by an ideal) the set of Cauchy sequences is given a purely algebraic description (again through category theoretical 'limits'). There is also the theory of associated graded rings which expresses the tangent cone at a singular point of an algebraic variety.  Artinian local local rings are the simplest of all rings, they correspond to points (the Hegelian One), or rather, infinitesimal points which yet have structure.

The theory of meadows, which arose in the study of abstract data types expresses the passage of ring theory (as part of Being) to a part of Essence. A meadow expresses the community and sharing of different rings united by a single point of infinity (analogous to compactification) - this is a very elegant expression of the categories of existence and substance.

In a Peircean terms we can associate Being with monadic predicates, Essence with Relations  (or diagrams in a category) and Concept with ternary predicates or indeed with higher order logic (maybe higher category theory), for instance involving the relations between an object and a relation.

\section{The transition to Essence and Modern Physics}

Generally speaking modern physics has a very Hegelian structure.
 In Newtonian mechanics we can start with simple postulates about invariance under Galilean transformations and the conservation of angular momentum 'emerges' as a consequence of the invariance of inertial frames under rotation. Starting from the invariance of the speed of light in special relativity we get that energy and momentum are united in a elegant way as a single 4-vector (which thus 'emerges') and that special relativity subsumes classical mechanics as a limit case when $c \rightarrow + \infty$ (Aufhebung).
There is also a further process of 'abrogation' in which concepts first considered static and absolute are then revealed to be relative. Being there (Dasein) is only being there through negation of the 'other'. This is exactly what relativity does by making time, length and mass no longer in-themselves but for-the-other (the observer). In classical quantum mechanics, starting from basic postulates about the wave function, observables and the evolution of the system we obtain that by passing to the classical limit we get an interpretation of the operators and the wave function in terms of classical Hamiltonian mechanics (including the wave function expressed in terms of the action $S$). This is another case of  'Aufhebung'. The uncertainty principle (a surprising result) is also case of emergence. The uncertainty principle interacts with the classical concept of angular momentum to give rise to an entirely novel concept / category: that of spin. Quantum mechanics (and specially quantum field theory) often abrogates previous concepts. Even in Bohm theory the concept of an individual particle is no longer tenable once we consider relativistic versions of quantum theory.

The narrative about the transition between a Newtonian absolute frame of reference and  the relativity of referential frames is rich in philosophical consequences.  Essence is negativity (this can be represented by the functor that collapses a category into the single category), the abrogation and vanishing of being which paradoxically is the retreat of being into itself, being reflected in itself, shining in itself, the realm of relation and mediation and positness.  The Newtonian absolute frame of reference  corresponds to the naive immediacy of being.  Relativity negates this: this frame of reference does not exist, no frame of reference can be an absolute immediate representation of physics.   Thus this negation carries over to all frames of reference, for each one of them relativity says: this is not physics.  Physics must stand beyond all frames of reference, the depth beneath the surface appearance. As Gödel defended,  relativity entails a certain form of phenomenism or idealism.  However the essence must manifest,  and from the negation of the absolute claims of reference frames follows that the essence must lie in the relation between the reference frames, that is, the movement and relativity between moments of Being,  And in physics this web of relations corresponds exactly to representations of Lie groups (or principal bundles).  And we have the further duality between the group itself and its representations which correspond to positing, Appearance, manifestation, Existence.  The sections about Identity and Difference at the beginning of Essence are all about relations and specially equivalence relations.   Physics can be considered as defined by equivalence classes modulo all possible frames of reference.

\section{The Miracle of Double Negation and Essence}

The beginning of the section of Essence in the Science of Logic is acknowledged to be difficult.  While the shorted version of the Encyclopedia Logic can be seen more down-to-earth, more connected to concrete embodiments of Essence as pertains
to the relativity of being, to classical scientific and metaphysical notions (maybe somewhat in the vein of Whitehead's Process Philosophy),  the treatment of the beginning of Essence in the Science of Logic (specially the section on Schein and the determinations of Reflection)  is of a very abstract nature, at once illuminating and frustrating, but suggesting many interesting connections to Hume, Pyrrhonism (skepticism is mentioned in the first remark as being superficial with regards to appearance) and Madhyamaka philosophy.  It is almost as if Hegel is struggling to articulate  a kind of  all-encompassing vision of the impermanence, relativity, lack of svabhava, of the world - all of which leads to the idea of something reflected into itself, shining into itself. The section on Schein can be construed as a metaphysical meditation  on change (but also on dialectic itself). It is also Hegel's solution to Fichte's problem of the independence of the laws of thought. Hegel deduced them all and expounds how they are all co-implicit thanks to his doctrine of absolute negativity and reflection.

The basic idea is that of pure or absolute negativity which can only maintain itself  through self-negation. Negativity posits itself through its own self-negation, that is, negativity creates a semblance, a shining, a reflection of being (at once presupposed, for there has to be something to be negated, and posited, for it is also a result)  in the process of its own self-negation.  We recommend the reading of the classical books by Hibbens, Harris, Noel,  Burbidge and Carlson and D.F. Ferrer's  very detailed discussion of Essence and Reflection in his book A Lógica e Realidade em Hegel (not translated into English).

Thus at the heart of Essence is the negation of negation, double negation.  But double negation is one of the richest and deepest of fundamental concepts in mathematical logic, theoretical computer science and category theory.   In type theory and computer science it turns out to express the idea of control, of a process having memory and its ability to travel through time and restore a previous state: this is at the heart of the doctrine of Essence.  Double negation is tied to the intrusion of imperative programming, and surely the idea of command is linked to that of negation. All command expresse negativity (do this ! means equally: do not do that !) and depends of negation.

We now give a brief description of how double negation arises in computer science connected to the notion of control and jumping through time and memory.

From a type-theoretic point of view there is no one single 'classical logic'. Rather there are many different classical-like extensions of minimal (or intuitionistic) logic, each corresponding to adding different classically but not intuitionistically valid axioms (DNE, LEM, Peirce's law, etc.). But we do not want mere axioms because they have no computational content and cause problems with term normalization.

Historically some functional programming languages started introducing control operators (typical of imperative programs) which can help handle exceptions or optimize programs. The idea is that in evaluating a term when reaching a subterm involving this operator the process can conditionally jump out (abort the rest of the evaluation of the subterm) of the normal evaluation process and 'pretend' that the subterm had some other value.

The structure of control operators is related to a style of programming called continuation passing style in which certain functions (which correspond to continuations) are always passed as parameters.

It turns out that if we try to type control operators (such as call/cc in Scheme, or try-catch, or abort, etc.) then these inhabit non-intuitionistically provable types.

And it turns out that writing a term in CPS is equivalent to the Kolmogorov double-negative translation of CL (as a type theory with control) into IPC.

Just as a fragment of IPC is isomorphic to CL, so too CPS style terms can be seen as programs with control (because continuators (which are like the state of the stack) are always passed as arguments)\footnote{for an introduction to this subject we recommend F. Agestam's Interpretations of Classical Logic using $\lambda$-calculus.}.

\section{Essence as Systems Theory and Holology}

 We should seek an abstract treatment of the category of Essence (Wesen) whereby Being negates itself to enter more deeply into itself and prepare the way for becoming in and for itself in Concept. The categories of Essence correspond to the categories of classical metaphysics, but they also corresponds to Objective Spirit, to nature. In modern terms the proper abstract treatment of the category of Essence (if it is not to lapse back into Being) must involve among other things a general systems theory adequate for theoretical biology. For it is the living organism which gives us light on the categories of existence (as organism and system) and active substance. Living beings manifest, they have appearance and phenomena. They preserve themselves in their own negation. Biology cannot be reduced to physics and is essentially more than physics in the same way that mediated Essence cannot be reduced to immediate Being. Life is process and mediation.  The theory of computation on one side is connected to formal logic. But on the other side it is a crude approximation of biocomputation, the computation of nature. What is a suitable mathematical model for biological systems that cannot be reduced to the general scheme of mathematical physics (sheaves, gauge potentials, jet bundles, PDEs, etc.) ?   Biology cannot be reduced to physics because it cannot be described by sections of a bundle which evolve according to a closed manifold in a jet bundle, that is, exhibiting the local causality of PDEs or integral equations based on these assumptions. The 'space' involved may be such that the stalk at any 'point' (and the subsequent evolution thereof) depends on the information on the entire space. For instance $X$ can be seen as a groupoid and we can think of the groupoid algebra product which for two elements is defined in terms of the entire groupoid.  Let physics correspond to a temporal process $\phi: T \rightarrow \Gamma(S)$ where $\Gamma(S)$ are sections over some bundle-like structure over a 'space' $X$. This structure forms a certain kind of (commutative) algebraic structure. Then biology corresponds to a temporal process $\beta : T \rightarrow B$ in which $B$ is a (non-commutative) algebraic structure which is not isomorphic to sections of a bundle-like object over any 'space' $X$. However $B$ may have phenomenological 'representations' or 'projections' $\pi: B \rightarrow \Gamma(S_\pi)$ for $S_\pi$ a bundle-like structure over some space $X_\pi$. Thus living organisms have physical measurable manifestations but cannot be reduced to such.

 Perhaps $\infty$-groupoids might play a role.  Here we discuss some of the intuitions involved.  Consider points in the plane.   Internal equality between two (externally unequal) points $a$ and $b$, written $a=b$,  is seen as being witnessed by there being a continuous path from $a$ to $b$.  In this way reflexivity, symmetry and transitivity have natural interpretations (in fact we obtain a groupoid or more correctly a 1-groupoid). Note that there can be more than one path between two points.  For a space $C$ consisting of two disjoint 'chunks' $A$ and $B$ there is no path from a point of $A$ to a point of $B$.  There are two extreme situations. One is the discrete case in which no two externally distinct points can be connected by a path - this is the classical concept of a set. The other extreme is a connected (technically: 1-connected) space in which we can find a single point $p$ such that $p$ can be connected to any other point in the space.  Thus the discrete space has the minimum of internal unity and the contractible space the maximum of internal unity. In general we can have a space with a certain number of connected components ('chunks'). If we view internal equality as an equivalence relation (and thus the space as a setoid)  then the number of components is the number of equivalence classes.

We now explain how being connected is only the maximum of unity at the lowest level. Given two points $a$ and $b$ there is a space of all paths from $a$ to $b$ denoted by $Path(a,b)$. Only now  in such a space  'points' are paths and the 'paths' are deformations between paths (homotopies). So given two paths $p_1,p_2$ in $Path(a,b)$ we have that $p_1=p_2$ is witnessed by there being a 'path' (i.e. deformation/homotopy) between paths $p_1$ and $p_2$. Now the same story is repeated at a higher octave. $Path(a,b)$ itself can possess the extremes of being discrete (set-like) or connected or anything between.  From elementary algebraic topology we get examples of spaces which are connected but which for two points $a$ and $b$ there are paths $p_1$ and $p_2$ in $Path(a,b)$ which cannot be deformed into each other (for instance on a torus). In the case where $Path(a,b)$ is connected for every $a$ and $b$ the space is called contractible but we see now why we could also call it 2-connected. There are also closed paths or 'loops', elements in $P(a,a)$, and there is no a priori reason why these could not be discrete, not mutually deformable into another. Thus we might add something to our characterization of classical 'set'.  If we think of a set as a collection of points or atoms, then these should not have interesting internal structure or symmetries and thus we  add the condition that $Path(a,a)$ be path-connected for every $a$: that is, all loops are internally the same.  We can now continue this process considering deformations of deformations and so forth. If a deformation is a path between paths in $Path(a,b)$, i.e. an element in $Path(p_1,p_2)$ where $p_1,p_2: Path(a,b)$, then we can also consider for $d_1,d_2 : Path(p_1,p_2)$ elements of $Path(d_1,d_2)$, that is to say, \emph{deformations of deformations}.  If we can keep iterating the path construction, then we obtain what is called an $\infty$-groupoid (this is only a rough description) which has a specific notion of $n$-connectivity for each level $n$. What kind of $\infty$-groupoid will have the maximum of connectivity or unity, i.e. how can we define a contractible $\infty$-groupoid ? Clearly by requiring at each level that $Path(x,y)$ for any $x$ and $y$ be path-connected, that is, that we have $n$-connectivity at each level $n$.

Now let us see how $\infty$-groupoids  might illustrate the metaphysical categories in Essence as well as the philosophy of nature in Hegel.   There is a lowest level of inorganic matter which corresponds to classical sets, thus to  0-discreteness. Then living beings have an organic unity which corresponds to 1-connectedness. Each part of a living body is connected to every other one on the lowest level, but can have quite distinct modes of connection.  We can also think of this as there being the same life present in every part but functionally  specialized. Next we come to mind which is a higher-order organism in that the same mind is present in every part in an essentially similar way (perhaps related to metacognition or "self-luminosity"), thus we would postulate  both 1- and 2-connectedness for instance. The failure of 3-connectivity could reflect the presence of time or sensation. We could also follow ancient metaphysics such as Plotinus and characterize the \emph{nous} (seen as a unified space of all forms) as being $n$-connected for all $n$, thus a contractible $\infty$-groupoid, the ontological maximum of unity.

A fuller understanding and interpretation of the section on Essence and its relation to modern mathematical physics will be attained through thorough analysis of our use of differential models in science as well as the allied concepts of determinism and computability, specially with regards to the possibility of extending these concepts to other more general mathematical frameworks.

\section{The Dialectic of the Combination of Systems}

Given two axiomatic-deductive systems $S_1$ and $S_2$, each is an island unto itself, each exists in the world of its own rules, the symbols of which are meaningful only relative to these rules and to the system as a whole. Thus $S_1$ and $S_2$ are neither mutually consistent nor non-consistent - and we can combine them into a new system in a free way.  To speak of consistency we must first identify symbols of both systems, we must somehow postulate a common language.  And we need a notion of contradiction. But all this only makes sense relative to a third system S3 which acts as an arbiter.  If $S_1$ and $S_2$ deduce contradictory things regarding a symbol $A$, how can we formalize this in $S_3$? 

By introducing a symbol for $S_1$ and $S_2$ inside the combined system $S_1+ S_2$, qualifying symbols representing 'from the perspective of $S_1$' and 'from the perspective of $S_2$'.  The introduction of subscripts or symbols for indexing by $S_1$ and $S_2$ is a projection of the system inside itself, a reflection-into-self.

But then $S_1$ and $S_2$ must be qualified by the arbiter system $S_3$. Why is the perspective of $S_1$ relating to $A$ what it is ? (And the same for $S_2$).  $S_1$ and $S_2$ share context $S_3$ and the quality of $S_1$ depends on that of $S_2$ and vice-versa.  $S_1$ and $S_2$ have become parts of a whole, for-themselves so far as they are for-another. Communication is the outcome of contradiction.

\section{Essence and Metaphysics}

Essence becoming actuality in the Science of Logic: beyond different models for a theory, different frames of reference for physics, representations for an algebraic structure, etc. Do we not here have the phenomenological reflection which neither attempts to have empty thought grasp directly its own structure nor is thought loosing itself in the objectified engagement in its action, but rather is the self-reflected awareness of thought in its thinking, a shift of perspective which knows itself in its process?  This goes beyond causality, computation and formal logic to inner and infinitel knowledge. Knowledge that what is in consciousness is taken and proceeds from its inmost center, a revelation, veiling-unveiling, which must pass to reveal the process: the revelation is re-velation.

Hegel's treatment of modality in the section of Actuality is very relevant to contemporary discussions on possible worlds and arguments against the world being a sum of independent contingent facts. There are also connections that could be made to modal logic (i.e. the connection between possibility, actuality and necessity and the axiom $\lozenge P \rightarrow \square \lozenge P$). 

The considerations on necessity and substance (and indeed the long digressions on Newtonian forces already suggest this) can be interpreted in terms of the solutions of differential equations (or in general smooth vector fields on a manifold). The local existence and uniqueness of a solution for a given initial conditions expresses determinism which (ideally, though not in mathematical reality) extends to a global determinism. 

The first part of Schopenhauer's World as Will and Representation is divided into four books. It would be interesting to establish a correspondence with the Science of Logic.  Clearly books 1 and 2 correspond to Being and Essence while books 3 and 4 correspond to Concept. However books 1 and 2 are written from the perspective 'for us', from the point of view of transcendental reflection. From this perspective books 1 and 2 belong to Essence and their philosophical knowledge expresses with great rigor and detail Hegel's theory of appearance, existence, substance, necessity, causality and actuality.  Or maybe Hegelian substance is much like the spontaneous production and manifestation of the will in Schopenhauer. It goes beyond phenomena and appearance because the will is immanent in its manifestation, totally in each one and beyond any particular one. The will has to manifest. Transcendental reflection, artistic contemplation and spiritual development represent the will's progressive self-knowledge, return to self, the Hegelian Spirit in the form of the concept.

There appear to be very interesting mutually illuminating connections between Schopenhauer's very abstract, universal, apodeictic and relational formulation of the principle foufold root of sufficient reason (and its further development in The World of Will as Representation) and Hegel's Doctrine of Essence (specially the treatment of absolute relation, necessity and actuality). This is a topic worth exploring in future research.

Category Theory: a category expresses essence. Each object is a different mode whose determination is inseparable from its relations (morphisms) to the whole (all possible determinations). An object is an expression of the category and yet not the category (for there are other objects). At the same time this circumstance of the object not being the category is itself internalized and expressed as the object being itself the sum-total of its relations (morphisms) with all other objects (i.e. with the totality of the category). This expressed also that the object is the category. 

The Concept can  be seen as arising when we have the cyclic transcendental consciousness of metatheory. That is to say: we realize that transcendental presuppositions are in cyclic dependency: to do logic we need computation, to do computation we need arithmetic and combinatorics, to these we in turn we need computation and logic. The germ of this in the reciprocal causality at the end of Essence (which we have also associated with higher groupoids). And also, as we have discussed in more detail everywhere, the partial (meta)reflection of one formal system by another (or the formal system by itself) is clearly relevant here: it is the idea of the whole being implicit to some extent in the part. For instance in ZF set theory itself we can prove that there are non-equivalent models of ZF.

  \section{Difficulties in Approaching Concept (Notion)}
  
 Perhaps a key to understanding Hegel's concept of Notion and the concept of Individual, Particular and Universal and their relationship - is to see how these concepts are embodied in the movements of Spirit in the Phenomenology of Spirit. That is, we view Notion as describing the structure and dynamics of a community of consciousnesses united in different ways into a higher consciousness wherein each individual also is fully identified with this universal consciousness. Although in the first part of Notion Hegel uses the terminology of the formal logic of his day, what he is getting at is much deeper.  There is here the mirroring of the microcosm and the macrocosm, the human mind itself is like a community of minds. This suggests the theory of fractals as well as Alain Connes' non-commutative geometry - for these all concern interface between the local and global. There are two extremes: the global contained fully in the local (as in holomorphic functions) or the local completely emptied of any global information: this is the case of Penrose tilings which cannot be distinguished by considering only local information. Hegel in discussing Objective Notion mentions a system in which the whole is implicit in each particular determination: this is the case of holomorphic functions: the germ at each point completely determines the whole function (on a connected component containing that point) by analytic continutation and the whole completely determines the germ at each point. In general this is a case of a sheaf in which the restriction morphisms are isomorphisms. Holomorphic functions can be called holistic because they are uniquely determined by their germs, the whole is in the part. Change a holomorphic function locally and you change it globally. The global information on the function is contained implicitly in any local restriction or germ (which unfolds via analytic continuation). Thus we can say that the "same" germ is spread out across the whole domain of the function.
 The case of holomorphic functions corresponds to the global being contained in the particular (or local). The global is thus totally local. But there is also the opposite situation in which the local is emptied of any global information. The global is totally non-local. This is the case, as we mentioned above, for the space of Penrose tilings. It is impossible to distinguish between two tilings by looking at local patterns. This could also be brought into connection to non-locality in quantum theory.

  The transition between Essence and Notion, involves the culmination into systems of reciprocal causality.  It would be interesting to read this in parallel with Connes' discussion on how many classical constructions yield equivalence relations or groupoids for which classical analysis and geometry are inadequate - in which a conceptual leap is required by developing tools for a non-commutative concept of space (i.e. a more sophisticated generalization of space which upturns the classical relations between the local and the global). The natural topology on such non-commutative spaces defined as the quotient of an equivalence relation is the coarse topology. So the germ of each point is the space itself, in Hegel's terms, the Individual is the Universal.
  
  Fractals and fractal calculus may well be worth exploring as a formal hermeneutic tool for Hegel's Logic. In the Encyclopedia Logic Hegel speaks of Objective Notion as undergoing a division into parts each one which contains the totality.
  
   For physics and systems theory the most developed and relevant sections in the Logic are found in the section on Object (mechanism, chemism and teology) and the first parts of the section of Idea (with regards to process and life). For the relevance to topology and geometry this is a good point of entry, with special focus on self-similarity (maybe in quantum field theory ?) and goal-oriented systems. The final parts of Idea are clearly about self-reference, self-knowledge, self-reproduction, self-production.  Concept-Idea is both essentially self-referential (it is knowledge knowing itself, and knowledge that is being and being that is knowledge, and also process and goal)  and self-productive as well as generative, in the sense that it generates or emanates nature. 
   
  \section{Notion, Consciousness and Mathematics}

The something and the other. In the Science of Logic Hegel thinks of the idea of two things bearing the same relation to each other and indistinguishable in any other way. Now a good illustration of this situation in found in the two possible orientations of a vector space. Each bears the same relation to the other and there is absolutely no way to uniquely identify one of them in distinction to the other.  In the same place Hegel talks about the meaning of proper names and states that they have none.

Hegel's strategy. All of Hegel is based on the dynamics and structure of consciousness as it reveals itself to itself. And yet for Hegel this structure is at the same time 'objective':  mere  'subjectivity' is seen as a failure and partiality of objective consciousness in not living up to its full objectivity or actuality (or knowledge and realization thereof).

Thus ideality and infinity are the key structure of a being which has developed to a point of being a subject. 

Finitude and limit:  a closed set.  Yet the boundary also defines what the set is not (cf. our discussion on the Heyting algebra of open sets in our paper on the Topos in Aristotle). Thus the boundary contains implicitly its own negation.

Limitation : an open set (cf. analytic continuation). Limit outside itself.

The ought... manifestation of the germ of a space through a given open set representation of the equivalence class. Any particular one is insufficient and can be replaced by another which is also insufficient.

False infinity: taking simply the set (or diagram) of such  open set representatives.  True infinity taking the limit (equivalence class) or limit (in the categorical sense).

Degrees of interpenetration and transparency between the individual and the universal (and between self-consciousness and essence) in the Phenomenology of Spirit.  From the rudimentary form in the ethical life to: the universal is in the individual, the individual in the universal and the universal is in the relation between individuals and the relation between individuals is in the universal (mutual confession and forgiveness). 

Milne's notes on Class Field Theory has a quote of Chevalley (1940) which seems couched in Hegelian terms: the Abelian extensions of a number field are determined by the number field itself, the number field  'contains within itself the elements of its own transcendence'.

'Pure' category theory is very poor in results, like Hegel's poor, empty, abstract universal.  But once we consider the particular, like Abelian or Monoidal or Model or Presheaf categories, things become very rich.

\section{Mereology and Homotopy Type Theory}

To prepare the main discussion of this section we make some introductory remarks about  \emph{homotopy type theory} and \emph{higher category theory}.
Homotopy type theory, an extension of the constructive type theory of Martin-Löf, was first proposed by Voevodsky as an alternative foundations for mathematics suitable for computer-assisted formalization and verification. For each type there is a whole hierarchy different levels of "equality"  analogue to the series of homotopy groups $\pi_i(X)$ of a topological space\footnote{ The linearly orderered real interval $I$  plays a central role in homotopy theory. 
	The concept of connectivity associated
	to it is pathwise connectivity (which can also be extended to Topos Theory\cite{Wraith}). }.

Higher category theory is based on homotopy theory, thus on model categories but from the perspective of enriched category theory and simplicial sets. Instead of ordinary (i.e. set-enriched) model categories we work with model categories enriched over simplicial sets $sSet$ (seen itself as a category enriched over simplicial sets. When endowed with its classical (Quillen) model structure its fibrant-cofibrant objects are called $\infty$-groupoids).  Instead of presheaves of sets we work with simplicial set enriched functors from the opposite of a simplicial set enriched category to simplicial sets. A key result is Dugger's theorem, the simplicial presheaf analogue for combinatorial model categories of the characterization of sheaf toposes as left exact full subcategory localizations of presheaves (see the next section for a discussion of this) - itself a generalization of the presentation of an object by generators and relations.  So this could be interepted as topos theory done (homotopically) over $sSet$ rather than Set\footnote{ 
	It seems we can give a more geometric interpretation of the nerve of category $C$ (a canonical way of extracting an simplicial set ) which is defined usually in terms of composable sequences of arrows. For three composable arrows $f,g,h$ we think of $f$ and $g$ as being in the plane but $h$ directed perpendicularly into space. Then we get 3-simplex in $N(C)_3$ with faces $(f,g, g\circ f)$, $(g,h, h\circ g)$ and $(g\circ f, h, h\circ (g \circ f))$ where we view a pair of composable arrows together with their composition as a triangle, i.e. a 2-simplex. Our intuition is that the nerve of a category keeps track of all commutative diagrams and each such diagram is a geometric object}.

There are competing definitions of (models of) $\infty$-groupoids besides the $sSet$-based one which itself is only one possible choice for what are called \emph{shapes} (which include cubical and cellular sets). All this suggests philosophically that our different concepts and models of what a 'space' is are special embodiments of a single 'pure' concept which is yet to be determined.  Homotopy type theory connects $\infty$-groupoids to  types, spaces or 'sets' of proofs/functions/computations.

There are interesting connections that can be made between topos theory and  mereology\footnote{Can we express Aristotle's  basic logico-ontological suppositions in terms of modern mereology ? Let us take a generalized organicist approach to the existence of wholes.  We assume strictly finite models. Let us write $PP'xy$ if there is no $z$ such that $PPxz \& PPzy$.  Then for Aristotle there are only the following two basic kinds of instances of $PP'$.  The first kind includes that of a certain kind of part of a living organism or being a part of a concept (for instance the parts of geometric concepts).  The second kind involves the relation between an individual and its immediate species and a species and its immediate superior genus. Abundant evidence can be adduced that Aristotle had a very organicist or at least geometric/topological approach to concepts, patterned somewhat after his biological theories. Genera and species are not conceived of in a strictly extensionalist manner, but do have associated extensions (modally conditioned). Thus for Aristotle  the only true 'wholes' are individuals and species/genera.  The status of parts of living organisms is an open question.  Also, are the only genuine parts of a genus its subordinate species ?  But surely the arc segments of a circle are not species of circle.   Our view is that this approach is inadequate because it conflates under the same concept both classical parthood and the species-genera relation. A more general framework is required with more than one mereological primitive. } (for a good reference see \cite{var}). Homotopy type theory also related to important problems regarding regarding identity, equality, intensionality and extensionality.

In mereology the parthood relation $P$ generally has the properties of a partial order which is a very special case of a category. Here we suggest that some core mereological conccepts make sense in the more general context of a  Grothendieck topology on a category.
Suppose we considered arrows $f :B \rightarrow A$ in a category as expressing generalized parthood and read them as $A$ is an $f$-part of $B$.   For $g :  C \rightarrow B$  then by composition $C$ is a $ f \circ g$-part of $A$.  First, let us investigate the concept of sieve, mereologically speaking.
A sieve on $x$ is a set (or class) of parts $y$ of $x$, $Pyz$, such that if $z$ is a part of $y$ then $z$ also belongs to this set.

\[Sieve(A,x) \equiv \forall y. ((y \in A \rightarrow Pyx) \& \forall z. Pzy \rightarrow Pzx) \]

A sieve on $x$ is thus a set of parts of $x$ plus all the parts of those parts.

A covering sieve is a generalization of fusion: $x$ is to be seen as a generalized fusion of any one of its covering sieves.  To justify this we argue that the three axioms for a Grothendieck topology (see \cite{Moerdijk} for the definition of Grothendieck topology) make sense for fusions. If $z$ is a fusion of the $\phi(x)$ and $Pyz$ then we can consider $\phi_y(x)$ (the $y$-restriction of $\phi(x)$) expressing elements which are overlaps of elements which satisfy $\phi(x)$ with $y$. Then we should have by the first axiom that $y$ is the fusion of the $\phi_y(x)$.

Suppose $z$ is a fusion of the $\phi(x)$. And consider a $P$-sieve on $z$ given by $\psi(y)$. If for every element $x$ which satisfies $\phi(x)$ we have by the second axiom that $x$ is a fusion of the $\psi_x(y)$ then $z$ is also a fusion of the $\psi(y)$.

Finally each $z$ is the fusion of all its parts: this is the third axiom

An important aspect of contemporary mereology  involves holology, the study of different types and degrees of wholeness in systems (we repeat here some of the discussion in section 12). The concept of higher groupoid is relevant to holology.   In Martin-Löf type theory  there is an extrinsic and an intrinsic concept of equality.  The intrinsic concept (written normally $\Gamma \vdash  t : A = A$) is what we shall mainly consider here. It can be given a topological (and higher-categorical) interpretation which is of some interest to mereology.  Here we only give a rough description of some of the intuitions involved.  Consider points in the plane.   Internal 'equality' between two (externally unequal) points $a$ and $b$, written $a=b$,  is seen as being witnessed by there being a continuous path from $a$ to $b$.  In this way reflexivity, symmetry and transitivity have natural interpretations (in fact we obtain a 'groupoid' or more correctly a 1-groupoid). Note that there can be more than one path between two points.  For a space $C$ consisting of two disjoint 'chunks' $A$ and $B$ there is no path from a point of $A$ to a point of $B$.  There are two extreme situations. One is the discrete case in which no two externally distinct points can be connected by a path - this is the classical concept of a set. The other extreme is a connected (or 1-connected) space in which we can find a single point $p$ such that $p$ can be connected to any other point in the space.  Thus the discrete space has the minimum of internal unity and the contractible space the maximum of internal unity. In general we can have a space with a certain number of connected components ('chunks'). If we view internal equality as an equivalence relation (and thus the space as a setoid)  then the number of components is the number of equivalence classes.

But being connected is only the maximum of unity at the lowest level. Given two points $a$ and $b$ there is a space of all paths from $a$ to $b$ denoted by $Path(a,b)$. Only now that 'points' are paths and the 'paths' are deformations between paths (homotopies). So given two paths $p_1,p_2$ in $Path(a,b)$ we have that $p_1=p_2$ is witnessed by there being a 'path' (i.e. deformation/homotopy) between $p_1$ and $p_2$. Now the same story is repeated at a higher octave. $Path(a,b)$ itself can possess the extremes of being discrete (set-like) or connected or anything between.  Indeed from elementary algebraic topology we get examples of spaces which are connected but which for two points $a$ and $b$ there are paths $p_1$ and $p_2$ in $Path(a,b)$ which cannot be deformed into each other. In the case where $Path(a,b)$ is connected for every $a$ and $b$ the space is called contractible but we now see why we could also call it, for example, 2-connected. We note that there are also closed paths or 'loops', elements in $P(a,a)$, and there is no a priori reason why these could not be discrete, not mutually deformable into another. Thus we could add to our characterization of classical 'set'.  If we think of a set as a collection of points or atoms, then these will not have interesting internal structure or symmetries and thus we could add to the definition of classical set the condition that $Path(a,a)$ be connected for every $a$: that is, all loops are internally the same.  We have seen the recurring motifs: connected and discrete. To this we can add external notions: there being no elements whatsoever in the space (type), there being exactly one according to external equality (singleton) and there being exactly two elements (for instance the type of classical Boolean values).

We can now continue this process considering deformations of deformations and so forth. If a deformation is a path between paths in $Path(a,b)$, i.e. an element in $Path(p_1,p_2)$ where $p_1,p_2: Path(a,b)$, then now we are given for $d_1,d_2 : Path(p_1,p_2)$ elements of $Path(d_1,d_2)$, that is to say, \emph{deformations of deformations}.  A structure like what we have here, in which we can keep iterating the path construction, is called an $\infty$-groupoid (this is only a rough description).  We get thus a notion of $n$-connectivity for each level $n$. What kind of $\infty$-groupoid will have the maximum of connectivity or unity, i.e. how can we define a contractible $\infty$-groupoid ? Clearly by requiring at each level that $Path(x,y)$ for any $x$ and $y$ be connected, that is, that we have $n$-connectivity at each level $n$.

A main idea of   the Univalent Foundations of Mathematics (also called homotopy type theory) is that mathematics is better done over $\infty$-groupoids rather than over sets. This brings unity between logic, category theory and geometry and also a natural way of developing computer software for formalizing and checking proofs. The concept of 'h-type'  is directly related to what we called $n$-connectivity.

The framework discussed above has applications to ontology and metaphysics - specially the Aristotelic and Neoplatonic systems, and in particular to the philosophy of biology and mind.  There is a lowest level of inorganic matter which corresponds to classical sets, thus to  0-discreteness. Then living beings have an organic unity which corresponds to 1-connectedness. Each part of a living body is connected to every other one on the lowest level, but can have quite distinct mode of connection.  We can also think of this as there being the same life present in every part but functionally  specialized. Next we come to mind which is a higher-order organism in that the same mind is present in every part in an essentially similar way (perhaps related to metacognition or self-luminosity), thus we would postulate  both 1- and 2-connectedness for instance. The failure of 3-connectivity could reflect the presence of time or sensation for instance. We could also follow ancient metaphysics such as Plotinus and characterize the \emph{nous} (seen as a unified space of all forms) as being $n$-connected for all $n$, thus a contractible $\infty$-groupoid, the ontological maximum of unity.

C. Ortiz Hill's book\cite{hill1} contains an interesting discussion on equality and identity in Frege and Husserl.  Per Martin-Löf and Frege played similar roles as founders of  radical new systems of mathematical logic and foundations of mathematics. As we just saw in Martin-Löf type theory we have two notions of equality. One is extrinsic and corresponds to (absolute) sameness (or identity) and satisfies Leibniz's law. The other is internal, and corresponds to equality, sameness according to some aspect\footnote{it would be interesting to compare this to Aristotle's statement of the law of non-contradiction in the \emph{Metaphysics}.}, which we interpret  in terms of paths.  Leibniz's law for internal equality becomes problematic and corresponds topologically to 'liftings' or 'transports'.

Extensional internal equality between functions is interpreted as a homotopy and we get the analogue of the notion of homotopy equivalence between two functions: each function is the inverse of the other up to homotopy (extensional internal equality). 

Husserl criticized Frege for conflating sameness and equality; extensional equality is a kind of equality and Law V expresses this conflation as an extensionality principle.
But Law V is powerful as a foundation of mathematics. We saw that Martin-Löf and Voevodsky agree with Husserl and keep the distinction between sameness and equality. However a certain controlled and transposed version of the Fregean approach becomes the Univalent foundations of mathematics. Instead of conflating sameness with extensional equality we conflate in a precise way internal equality with equivalence (which involves likewise extensional equality). More precisely the univalence axioms reads:

\[\text{\emph{Internal equality itself is equivalent to equivalence.}}\]

In the section of the determinations of reflection in the Science of Logic and the discussion of Identity and the Principle of Identity (which is critiqued) it is clear that Hegel is considering higher-order identity in that he considers the identity of identity itself and difference and so forth (which is also, be it noted, found in the discussion of the five supreme genera in Plato's Sophist).  Is seems that Hegel's concept of identity is implicitly not only higher order but homotopic. There is surely much to be explored as regards to the connection between homotopy type theory and the Logic of Essence.

\section{Aufhebung in Type Theory}

Consider the Barendregt Cube or more specifically the system of logics studied in Jacobs book 'Categorical Logic and Type Theory' which starts with simple type theory and culminates in higher-order dependent type theory which subsumes all the rest. There is a process in which each level subsumes the previous one in particular through reflection-into-self, the prototype being the subobject classifier wherein Prop becomes a Type. The culmination of the Barendregt cube is $\lambda C$. There is a parallel with Hegel's development of Concept (Begriff). In the simple typed lambda calculus we have the classical division between subject and predicate $ t : T$, $t$ is of type $T$, $t$ is a $T$, which corresponds to the division between term and type. The development along the directions of the cube expresses the increased dependency of terms on types, types on terms (and types on types) until the distinction is all but subsumed and abolished in $\lambda C$. There is only one pure 'type' $\square$ which expressed the Concept. Jacob's book allows us to trace this in terms of fibered categories.

We propose that we extend the Barendregt cube incorporating substructural logics via monoidal categories and generalisations of linear logic. This involves replacing the cartesian product with the monoidal product and losing weakening and contraction. Note that symmetric monoidal categories are a primary example of the internalization of the hom-sets of a category within that category.

Consider how predicate logics are extended via propositions-as-types. The poset (or preorder) for Prop becomes a full-blown category of proof-terms. This is like the Understanding in the Hegelian dialectic that wishes to erase the past, the temporal dynamic constitutive and genetic process of a concept. The understanding says what matters is what is proven not the process of the proof (repressed history). Against this dialectical reasoning forces the proof itself to be remembered and speculative reasoning incorporates the proof-process into a concept, thus turning propositions into types and proofs into terms.

\section{Kant and Sheaf Theory}

Lawvere's theory of smooth toposes and their use for formulating differential geometry and differential equations can be seen as an unveiling of fundamental synthetic a priori categories and principles. The adequacy and conformity to scientific experience (in the spirit of critical idealism) is a vital dimension in the philosophical deployment of category theory.

The sequence seems to be : topos $\rightarrow$ ringed topos $\rightarrow$ lined topos $\rightarrow$  smooth topos.  Our task is to elucidate the phenomenological and categorical (in the philosophical sense) meaning of the concept of topos.  Also to understand why the internal version of the concepts of commutative ring, commutative algebra, linear map (cf. Weil algebras and the Kock-Lawvere axiom) and an infinitesimal version of simplicial objects (used for defining the dg-algebra of differential forms) appear to be of so fundamental and vast a scope as a priori conditions for mathematical physics and other branches of natural science. And how are we to understand Kant's concept of space in the light of synthetic differential geometry, it being  point-free ? Seeing a vector field on $X$ as an infinitesimal deformation of the identity map on $X$ seems very intuitive.

Kant's theory of schematism could be interpreted in particular as implying that any line-shape $R \rightarrow X$ must be seen as a solution of a differential equation; we give the Weil algebra a generative dynamic interpretation.

Let us look at the first two antinomies of pure reason in the transcendental dialectic, involving time and composition. The point-free nature of space is presupposes for the argument to work. Can causality involving a previous moment of time be captured in terms of the infinitesimal path $D$ via a prolongation principle ?  The infinity which cannot be object of a synthesis appears to be best captured by a set dense set for $<$.   Otherwise why cannot for instance the points $\frac{1}{n}$ be objects of a synthesis ?

Kant appears to be saying that every state of the universe must have a temporally previous state but at the same time there cannot be infinitely many previous states to a given state.  

We must compare the antinomies to various intuitionistic principles and classical principles not intuitionistically valid.

It could be tempting to see a sheaf as a derived concept relative to an abstract bundle (which is a very Kantian notion). Sheaves are sections of some bundle. But we must check if this carries over for Grothendieck topologies too, that there is an equivalence of categories between sheaves over a site and étale bundles defined not in terms of ordinary topology but in terms of Grothendieck topologies (as in étale cohomology).  However there are other entities in geometry which are more naturally seen as sheaves than as a bundles: for instance the sheaf of continuous or holomorphic functions. A sheaf is essentially a phase space, a space of phenomenal possibilities which expresses how these possibilities flow locally and cohere. But in actual situations the number of different sheaves (over a given space) is very definite and determinate. For instance, sheaves of smooth functions on a manifold and smooth $k$-forms and other functions relevant to analysis. Many of these actual sheaves form natural complexes, so much so, that the complex itself can be seen as adequate kind of generalized sheaf (cf. the derived category).  Since it is the cohomology of the sequence that mostly of interest, complexes are identified if they have the same cohomology : this is the basis of the derived category constructions which turns morphisms of complexes module homotopy into 'fractions' where  in the denominator quasi-isomorphisms are inverted.

But how are actual sheaves, sheaves of concrete interest, 'generated' , beyond the basic ones discussed above ? By functors generated by continuous maps $f: X \rightarrow Y$.  Of great interest is the study of neighbourhoods of fibers $f^{-1}(x)$ as $x$ varies.  $Y$ is often seen as parametrization space or base space. Hence the presheaf  on $Y$ given bt $V \rightarrow \Gamma(f^{-1}(V), F)$ for $F$ a sheaf on $X$. This carries over to a derived functor $Rf_{\star}$ taking complexes of sheaves on $X$ to complexes of sheaves on $Y$ by which is studied the cohomological variation of the fibers along $Y$. This is a main source of the generation of interesting sheaves, the sheaves used in practice.  Another source of is the functor $f^{-1}$ by which is studied the cohomology in an infinitesimal neighbourhood of the images of open sets via $f$. These two functors are abstract versions of integration and differentiation. Another important operation is the "restriction" composition $R_\star j j^{-1}$ for an inclusion $j : V \rightarrow X$. This takes a sheaf $F$ and yields a new sheaf which, roughly speaking, takes into account only the nature of $F$ on $V$ or in infinitesimal neighbourhoods of $V$.

Thus a category of sheaves becomes interesting and intelligible by its relation to other categories of sheaves.  So considering at once the 2-category of sheaf-categories, or more generally, toposes, is very natural and imposes itself in the nature of things\footnote{The adjunction between toposes and locally presentable categories is discussed in Marta Bunge, Aurelio Carboni, The symmetric topos, Journal of Pure and Applied Algebra 105:233-249, (1995).}.

Having a whole category of sheaves leaves a vast amount of elbow-room. A category of sheaves represents a spectrum of different spaces of possible manifestation (cf. how locally constant sheaves can be identified with covering spaces). The category of sheaves over a given topological space represents every mode of phenomenal possibility space of that space - thus in a way the category can be identified with the space itself.

A remarkable property of sheaves is their homogeneity for scaling. Given an open set $U \subset X$ we get automatically from a category of sheaves on $X$ a new category of sheaves on $O$, $\Gamma_U : Shv(X) \rightarrow Shv(U)$.

Categories do not apparently have the vertical hierarchical structure of the classical genus-species classification. For instance: group is a species of monoid and abelian group is a species of group. There are corresponding categories Mon, Grp and Ab which form a chain of subcategories. Our construction from our paper can be interpreted in terms of successively taking equivalence classes of equivalence classes of equivalence classes.  In category theoretic terms this translates as a sequence of categories $C_0, C_1,C_2,...,C_n$ and a sequence  $F_1, F_2,...,F_n$ of sets of morphisms in $C_{i-1}$ such that $C_i$ is the quotient category of $C_{i-1}$ via $F_i$.  Consider how a given infima species might be described in terms of a protoype $P$, membership  of an object $X$ being ensured by the existence of a deformation $f : P \rightarrow X$. We can also think of an interpretation in terms of $\infty$-groupoids but it is more subtle; it is a top-down approach using connected components and $n$-contractibility.

But sections of a sheaf are like individuals of a given infima species. The category of sheaves is like the genus of the infima species. Then the 2-category of (small) sheaf categories is like the yet higher genus of this genus. Functor categories are like the category of relations.

So: section $\rightarrow$ sheaf $\rightarrow$ sheaf category $\rightarrow$ 2-category of sheaf categories. Aquinas (in De Ente et Essentia) views genus as a space of possibilities rather than as a minimal matter to which difference is added as a form. Difference determines or 'picks out' a latent possibility of the genus. This agrees with the concept of sheaf and category: we choose a section of a sheaf or a definite sheaf in the category of sheaves.

Formal axiomatic metaphysics in the sense of Zalta all depends on a multi-valued logic having truth-values $\Omega^W$ were $W$ represents possible worlds.  Or rather, consider that $W$ must be endowed with a topology or Grothendieck topology so we can consider $\Omega$ for $Shv(W)$.

Thus we have a fundamental axiom of metaphysics: the set of possible worlds must be endowed with a topology. The idea of alternative situation always depends on a more or less strong continuous deformation of the actual world. If a sentence $P$ holds at a world $w$ then it must hold in some neighbourhood $U(w)$ of $w$. Thus it is natural to capture modal logic by the topos of sheaves over $W$, $Shv(W)$.

In the section on understanding in the phenomenology of spirit there are many interesting considerations about 'force' and exteriorization and the unity of motion determined by a law.  The germ of a sheaf is the unity of between interiority and exterior manifestation (for there is not definition for the point itself independently of its neighbourhood).   Let $W$ be a category $\mathbb{W}$ with a Grothendieck topology.  Then take modal logic to be a functor $ m: \mathbb{W} \rightarrow \mathbf{Cat}$.

Topos theory is (very roughly) the study of how logic coheres and varies with space and time and  possibility. All predicates $X\rightarrow \Omega$ have a domain which is a 'type' $X$.  But in general $X$ is a sheaf. Types as spaces can be expressed as 'types as sheaves'.  For instance a type parametrized by possible worlds. Thus predicates are in this case coherently parametrized according to possible worlds. But notice how in general the objects $X(U)$ can vary for $U$ in the base category. If $X$ were an atomic type of sets of individuals then the actual individuals could very according to the possible world. An equivalent way to see this is as the subobject fibration in which we view propositions on $X$ in correspondence with Sub(X).

Zalta's encoding could be descriped as a morphism $enc_X : PX \rightarrow PX$.

Jean Hyppolite lays emphasis on Hegel's positing of  'the identity of identity and difference'. But it is difficult not to think here of the univalence axiom: 'the equivalence of equivalence and identity'. Hegel's logic, despite this going against the surface intention of Hegel himself, may well be capable of a formal axiomatic treatment. This will depend on a proper account of identity and equality.

In Jacobs' Categorical Logic and Type Theory there is the idea of giving a categorical semantics for untyped lambda calculus ($\lambda$-categories on p.155-156) related to Scott's reflexive objects. We take an object $\Omega$ in a Cartesian closed category for which $\Omega = \Omega \rightarrow \Omega$. This expresses that $\Omega$ has a mediation within itself, is self-mediating. Is in and for itself.

If conscious experience is normally present in unreflected 'globs' , the goal of analytic insight is to unmask and be continuously aware of the ternary structure present in consciousness $\bullet \rightarrow \bullet$ and its subsequent higher-order unfoldings.

We mentioned before the archetypal structures and processes of consciousness. Here is an incomplete tentative list (with an implicitly Kantian basis):

Synthesis - gluing, covers, the sheaf-condition = extensibility on $j$-dense objects for a topos with a Joyal-Tierney topology.

- different orders of wholes (higher groupoids)

Self-reflection - a system which can represent (partially at least) higher order aspects of itself within itself.  This is the original synthetic unity of apperception = I know that I am knowing. This is found in recursive definitions, inductive types, the successive powers of the $\lambda$-cube wherein external aspects of the system become internalized and internally represented, also the subobject classifier, truth-value object $\Omega$ in a topos. 
Return-to-self, that is, Kant's trinary structure in the CPR.  This is related to the negation of the negation, double negation as the third (synthesis).  In topos theory this relates to the dense topology and in particular to forcing.  The idea is simple. In rough terms it is as follows: consider $U\Vdash\phi$ as signifying that the sentence $\phi$ holds in region $U$. We define $U \Vdash \neg \phi$ if $\phi$ does not hold on any subregion $V\subset U$.  Then $U\Vdash \neg\neg \phi$ means that for any subregion $V\subset U$ we choose we must have that there exists a $W\subset V$ such that $W\Vdash \phi$.

Double-negation can also be connected to temporality: something must pass to reveal itself, ti to einai, quod quid erat.

\section{A Recapitulation}

Let $C$ be a category representing the totality of individual consciousness and their relations.  Then for a given individual consciousness A we have  i) the set $hom(A,A)$  (individual, internal relations),  ii) the set $hom(A,B)$ relations between two individual consciousnesses, and finally iii) the slice categories ($C \downarrow  A$) and ($C \uparrow A$) whose objects are all possible relationships between A and its environment (all the other objects in the category $C$) and whose morphisms are relationships between objects "from the point of view of $A$". Also the process could be repeated, taking slice-categories of slice-categories and so forth.

Notice how holistic the definition of limit and colimit are.  The property of a cone $C$ being the limit of a diagram depends in an essential way on the entire enveloping category (i.e. a universal property)

Elementary toposes represent consciousness and the mind while Grothendieck toposes (toposes of generalized sheaves) are a very abstract representation of physics (and perhaps also systems theory and biology).  The  underlying site (or in particular topological space or groupoid) represents the most general concept of space-time, the sheaves on it are generalizations of the fields in physics.

How are we to understand from a Hegelian perspective two objects being isomorphic ? It is like the self and the other become mutually determining in the 'true infinite' , or the reciprocal causality before the transition to concept.  A connected groupoid in which all objects are isomorphic to each is much like 'quantity' in the Logic.   In a quotient category we introduce new morphisms thereby forcing certain morphism to become isomorphisms.  Taking the quotient is like moving forward in the processing through the categories in the Logic.

Now the section on Measure on the Logic is highly suggestive from the point of view of modern topology. For modern differential topology and singularity theory represents a development of  key intuitions and concepts found in this section.

First of all: the relationship between quality and quantity recalls that not only of physical field spread out through space-time, but the general concept of a sheaf. The base category is often space, or space-time, or some control space, while the target category is phenomenological, qualitative although having of course having in turn a quantitative aspect.  The base is extensive, the fiber intensive. The further development in this section of the Logic - and most specially the concept of Nodal Line - corresponds very closely to Thom's general theory expounded in Mathematical Models of Morphogenesis (which also has close connections to Aristotle). There is a strong connection here to my paper on continuity and topos in Aristotle (which was published in Intentio). On the other hand, Hegel offers an important insight on the essential co-dependence between the qualitative and quantitative which has consequences for modern topology, singularity theory, differential equations and dynamical systems.  Strictly determined systems have to be studied qualitatively. But qualitative studies are only meaningful and fruitful in the context of a specific determined object. In order words the topological exists to study the metric and the metric has to be studied using the topological. The same goes for differential equations and dynamical systems and singularities. As Hegel wrote: the quantum is quantity (the topological, in the indifferent) endowed with quality.

The discussion on the sorites and 'phase transitions' is prescient. It recalls  locally constant sheaves (covering spaces) - for instance if we analytically continue the complex square root around the origin we do not end up where we started, in terms of the section /quality...

Hegel also offers interesting considerations on the philosophy of physics (and philosophy of science in general). Is not general relativity a partial realization of Hegel's project of deriving physics purely from the properties of space and time ? Hegel's observation of the influence on the observer on the experiment (in his discussion of thermodynamics) is prophetic.  The transition to essence suggests that while applied mathematics is legitimate in its own domain and perspective, deeper scientific knowledge has to be found elsewhere - but the latter must pass through the thorough understanding of the former.

Hegel's concept of type is very interesting. It suggests that it is like the image of a continuous map $F: U \rightarrow \Phi$ where $U$ represents the allowed variations of the independent parameter. $F(W)$ a type much like in homotopy type theory, because endowed with a topology and a concept of homotopy and continuous variation. Besides the connection to Thom there is also the more general connection to moduli spaces.

We need to delve deeper into Hegel's concept of indifferent quantity. This is clearly the concept of germ around a point: a fundamental system of neighbourhood, a direct set, etc.

Open and dense.  This is related to tranversality - again nodes. And to stability, expresssion of mediated quantity and invariance of quality. We can also investigate   the topological aspects of the calculus of variations.

Consider an analogy : day is to night as life is to death.  Surely in any analogy there is implicit a correspondence, a set of functions which takes the pairs to each other: \\

F: {day,night} $\rightarrow$ {life,death}\\

 and\\

 G: {life,death} $\rightarrow$ {day,night}\\

  Consider the statement: death is the night of life. If this is so then surely:  night is the death of day and vice-versa. We take \emph{of life} and \emph{of day} to be the functors $F$ and $G$ respectively - here we venture that 'life' and 'day' designate the whole genus, the set of both elements of the respective pairs. Thus it would make sense to day : day is the life of day  and life is the day of life. And we read this as $death \rightarrow F\,night \Leftrightarrow G\, death \rightarrow night$ which (together with the case for day and life) corresponds to the definition of an adjunction (we can assume we are in a groupoid in which 'is' is an isomorphism). We can also write symbolically:

\[\frac{day}{night} = \frac{life}{death} \Leftrightarrow day\times death = life \times night\]

\section{Formalization of Hegel's Logic of Concept}

We now address the third great part of Hegel's Logic, the one dedicated to Concept (Begriff) and its unfolding into Idea. Here we give a few considerations on how the theory of Concept might be mirrored in formal logic, mathematics and science. Besides the direct interest to the interpretation of Hegel we believe that such investigations may be of interest to formal approaches to phenomenology and dialectics and shed light on the problem of formal models of consciousness.

One could object that Hegel states explicitly that all such formalization will pertain to the strictly finite exterior domain of the understanding and as thus be inadequate. We reply that modern developments present us with formalisms which embody multidimensional, dynamic, modal, meta-theoretical and self-referential aspects which were not evident in mainstream 19th-century logic and science.

The logic of concept deals with the phenomena of self-reference and the reflection of metatheoretical properties of a system within the system itself (including incompleteness).  Other key concepts will be self-similarity (Hegel in the Encyclopedia Logic writes of Object that it is a division into distinct beings each one of which is the totality),  process philosophy and skeptical and transcendental-critical dialectics.

But we should not mystify the section on Concept, nor in particular Subjective Concept.  It concerns a transcendental phenomenological analysis of very ordinary down-to-earth notions and cognitive structures - and this is a leitmotiv for Hegel: the most down-to-earth is given the most forbidding transcendental description.

So subjective concept is really about our plain old concepts and the logical practice of Hegel's day (in this case there seems to be a strong Kantian influence). And yet questions at stake are contemporary and of the utmost difficulty and profundity remains: What is a concept?  What is a judgment? What is inference? What is induction? What are natural kinds? What is the extension of a concept? Do proper names signify? What are definite descriptions?  What is an individual? Are there individual essences? How we to interpret the traditional theory of genera and species? 

Hegel's treatment of subjective notion suggests that the analysis of what concepts are is already a powerful clue to understand the formal structure of consciousness. 

We are in the presence of a transcendental reflection on concepts and the mind's process of forming, grasping and using concepts.  A concept divides itself up as a genus into species. But it also has an extension comprised of individuals. But individuals have their set of properties, concepts which inhere into them.  The individual in the comprehension of a concept can have other concepts inhering in it. Thus concept appears completely disjoint and divided into a "subjective" and "objective" component. The problem is uniting them and achieving an individual-universal .  The fundamental relation is the relation between the individual and a concept - like Zalta's encoding or Gödel's encoding or in general Lawvere's presentation of representation of concepts by constants. And in Hegel there is a subtle presence of the problem of the meaning of proper names, of individual terms and his position seems to be that of Mill.

We shall discuss ahead the importance Lawvere's work and of categorical logic in interpreting and formalizing Hegel's logic of concept. But we can remark for now that in the Encyclopedia Logic 191 regarding the syllogisms of necessity we can interpret the categorical syllogisms (in which the particular mediates) as a type judgment $t : A$ and $A:U$ for some universe $U$.  The hypothetical syllogism is linked by Hegel to the mediation of the individual and we can think of this, for instance, as some kind of term $\lambda (x:A) B: A\rightarrow B$: this functional term takes individuals of type $A$ and ones of type $B$ - it thus "mediates between individuals". The disjunctive syllogism is linked to the universal. We could interpret this in terms of $\Sigma (x:A) B(x)$ which represents the integration of the particulars $B(x)$ for $x:A$ into a single universal.

It seems to us that the comparison between the classical theory of definition involving genera an differences and modern logical and definitional practice is relevant to Hegel's arguments (see our paper "Modern Definition and Ancient Definition") as well as René Thom's geometrical-topological theory of semantic spaces and predication. Thom might have agreed that Hegel introduces ontological and subjective considerations into formal logic (cf. the essential distinctions between: the apple is red, the apple is edible, the rose is a plant, the rose is beautiful). This might suggest a geometric interpretation in which we distinguish for instance between accidental properties of a space (relating to fiber bundles and connections), relational properties such a cobordism and more intrinsic properties such a homotopy and cohomology.  Hegel's critique of formal logic recalls that of René Thom. 

We can also think of monads and commands in relationship to judgments. The unit of a monad integrates an individual into its particular (a wrapper) while the counit of a comonad extracts the individual from the particular (context).

In the treatment of judgment and syllogism there seems to be at play a dualism between comprehension and extension in the modern sense. 

How are we to understand the transition to object and the whole theory of objective concept? 

We can think of the universal as the center, the "self", which radiates intentions (particulars), directed rays proceeding out of the center which stop at points (individuals) expressing pure negativity, closed in on itself, inconceivable. We pass from a single ray from the center (or multiple rays to a single external point) to multiple rays from the center to multiple points and pass to the center with a series of looped rays which return to the center expressing the identity between the universal and the individual meditated by the particular - the genesis of object.

Also the following analogy occurs to us:  the universal is like a group $G$,  the particular is like a particular representation of $G$ on some vector space $V$ and  individuality is the vector space $V$.  The concept is the  monoidal category of all representations of $G$. For certain kinds of monoidal categories the group can be extracted from the category just as the category is generated from the group by considering its representations. This expresses the flowing forth between the universal, particular and individual.

But in terms of our formal interpretations what precisely is new compared to Essence? One could propose that subjective concept corresponds to syntax and formal logic while objective concept corresponds to models, to model theory.  In Hegel's treatment of mechanism, chemism and teleology there is certainly a more holistic, systems theoretic attitude involved together with a closer connection to formal logic.  Thus it is the section of objective concept that categorical systems (or process) theory (whether based on homotopy theory, monoidal categories or topoi and sheaves) will present itself as the most adequate and possibly illuminating formalization. It is here perhaps that groupoids and higher category theory would find their most adequate applications.   This patent in Hegel's text which is concerns precisely with the dialectic of a system of objects and their mutual relations and the problem of intrinsic nature being defined in terms of the their relations to the whole:  in objective concept the categorical theoretic elements in Hegel become manifest and explicit. Also the categorical approach allows to understand the connection to subjective (formal) logic through the logic of the category and in particular through the correspondence between theories and topoi.  The kinds of higher order logic (with partial elements) employed in topos theory (cf. the work of Fourman and Scott) are of particular interest though we can question (as will be discussed further ahead) the extensionalist, structural, finitary and non-paraconsistent assumptions.

If objective concept is certainly the highest place for a categorical general systems theory in Hegel we could also ask if concurrency and the $\pi$-calculus could play a role here. Can we formalize goal-seeking in a categorical systems theory? We could examine topological accounts of teleology given by the concept of attractor (or even organizing center) in dynamical systems theory.

Monoidal categories are certainly a promising candidate to express the integration and correspondence (the cycle) between the subjective concept (logic) and objective concept (general systems theory, process algebra, physics, geometry, computation).

We can also justify that the logic of concept is the category theoretically fully explicit and adequate part of the Hegel's logic. Consider a category $C$. Then the consideration of the objects as objects corresponds to the Logic of Being wherein morphisms connecting objects appear a something both necessary and external. The Logic of Essence corresponds to considering morphisms as morphisms wherein their source and target objects are implicit.  The identity morphism $id_A$ corresponds to the principle of identity, etc.  The transitions might be understood in a 2-categorical context. In more advanced stages one considers diagrams into $C$. Finally the Logic of Concept corresponds to considering the whole category $C$.  The logical aspect is the internal logic of the category $C$. In more advanced stages we can consider the presheaf category or sheaves over a site defined with $C$ (or an embedding into a cartesian closed category) and finally come to the category of categories or other considerations from homotopy theory. Maybe the particular can be seen to correspond to slice categories.

We propose that the final synthesis between subject and objective concept in Idea must be related to the overcoming of the classical duality between syntax and semantics. We find this in categorical logic (and in algebraic logic). A topos is not merely a "model" of a "theory". It contains logic and theory within itself (cf. the syntactic categories). Note that in the section on Life Hegel assimilates syllogism and process. This is an anticipation of the Curry Howard correspondence. A syllogism is a proof, hence a term, hence a program, hence a process. Cf. 243 "It thus appears that the method is not an extraneous form, but the soul and notion of the content". We could say in the same way that in intuitionistic logic the soul of a formula is not its truth or falsity but its proofs.

Reading the last sections of the Encyclopedia Logic we are lead to propose that the ultimate expression of the Hegelian Concept, the Idea,  can only be a formal system which can reflect its own metatheoretical properties (such as incompleteness)  as well as of weaker systems and which can reflect its own proof theory, it can internalize the very process by which the external subject acquires metatheoretical knowledge in the system: the great example being due to Gödel or in general the classical semantic paradoxes of self-reference and truth. "This sentence is false" can be given a Hegelian interpretation. On the immediate level this expression is a subjective, spontaneous, free, immediate positing "I am freely positing that what I am saying is false" and thus in the aspect of its immediacy and positedness it is true, because it only expresses a free choice and finite determination as such (cf. omnis determinatio negatio). However by reflection and return-to-self this determination dissolves itself, its determination, its truth-claim is the seed of its own overcoming.

Can type theory and category theory encompass this kind of reflection?  Or do we need to break free from type restrictions and alter some fundamental logical assumptions to be able to encompass things like Girard's Paradox (or the Burali-Forti paradox)?  Saying that $\vdash U: U$ is an expression of the Hegelian Idea.  $U$ remains the same in all its particularizations. 

The key to this will be Lawvere's unified treatment of diagonalization arguments, proof theoretic and aletheic incompleteness results and self-reference.  Lawvere's categorical treatment of metatheoretic representability is a high-level reflection on cognition itself, the key part of the final section on Idea.

Thus we propose a careful study of Lawvere's  1969 paper, Diagonal Arguments and Cartesian Closed Categories.  Can Lawvere's results be extended from cartesian closed categories to certain kinds of monoidal categories (compact closed)?  Can incompleteness be used to interpret quantum theory?  

Because Lawvere's categorical approach strips Gödel's theorem of its seemingly mysterious, purely syntactical appearance, he viewed the theorem as a beautiful demonstration of objective dialectical movement. Rather than a limit on what we can know, he framed it as a natural mathematical manifestation of the contradiction between an "entity's internal capability to represent" and its "external limitations".

McTaggart's commentary of Hegel's logic is of great interest even if frequently overtly critical. The same goes for Burbidge's earlier work on Hegel's Logic.  MacTaggart's interpretation of judgment is interesting in the light of the above work by Lawvere and Zalta's theory of encoding. McTaggart writes: the perspective that starts from the individual breaks down. We might interpret some arguments in the sections on judgment and syllogism as arguing against logicism, arguing that purely logical rules could never be the right kind of mediation:  there are no purely logical rules which assign properties to a contingent individual. Regarding 190, in modern type-theoretic notation the rule that from $\Pi(x:A)B(x)$ and $a:A$ we can deduce $B(a)$ is rejected as a syllogism due to both a distributive interpretation of $\Pi$ and the rejection of "tautologies" like $A \vdash A$. To Hegel the meaning of $\Pi(x:A)B(x)$ must coincide with the enumeration of all instances $B(c)$ for $c:A$ which will include in particular $B(a)$. The discussions related to induction seem relevant to the modern intuitionist logic with partial domains.

The Hegelian Idea is a kind of dynamic synthesis between Sextus' Pyrrhonian dialectics and Kantian transcendental criticism - but perhaps even closer to neoplatonic dialectics. Dialectics as the life of the nous which unfolds, divides itself, circulates through itself then reunites and integrates itself within itself. In the words of Hegel himself:\\
\

215. \emph{The Idea is essentially a process, because its identity is the absolute and free identity of the notion, only in so far as it is absolute negativity and for that reason dialectical. It is the round of movement, in which the notion, in the capacity of universality which is individuality, gives itself the character of objectivity and of the antithesis thereto; and this externality which has the notion for its substance, finds its way back to subjectivity through its immanent dialectic. }\\

What Sextus lacked was a metalogic, a methodological unfolding of his path of equipollence as it traverses the regions of logic, physics and mathematics.

Finally we mention that paraconsistent and infinitary logic might be of interest in the formalization of the logic of concept. Paraconsistency is suggested by Fichte's Wissenschaftslehre: the ego posits the non-ego and we have $A \wedge \sim A$. From a Hegelian and dialetheistic perspective it seems reasonable to opt for the presentation of paraconsistent logic which discards the negative syllogism $A\vee B, \neg A \vdash B$ as disjunction introduction (or the disjunctive syllogism) seems to be given a more favorable view. Infinitary logic would appear to express naturally the integration between the subjective and objective.

\section*{Appendix 1: Homotopy Theory and Higher Categories}

Note that we are not implying that category theory is preferable to set theory as a foundation for mathematics or even for science in general.  We are inclined to hold that dependent type theory is a good candidate for such a role.

From a categorical point of view what is the simplest object we can conceive ? The singleton category with only one object $\star$ and only one arrow, the identity morphism $id_{\star}$ on this object. All such singleton categories are equivalent and a singleton category is in fact the terminal object in the $2$-category $Cat$ of small categories. They are complementary because they emerge from right and left adjoints respectively (which are necessarily fully faithful) of the same unique functor $\mathbb{T}: \mathcal{C}\rightarrow \{\star\}$.

\[
\begin{tikzcd}
	\mathcal{C} \arrow[r, "\mathbb{T}"]  & \{\star\}\arrow[l, hook', "1", shift left = 2ex] \arrow[l, hook', "0"', shift right = 2ex]
\end{tikzcd}
\]

This gives rise to adjoint modalities $ \emptyset \vdash \star : \mathcal{C}\rightarrow \mathcal{C}$. We can see the singleton category as the most rudimentary being, the Etwas, something.

The first section of Hegel's logic, the logic of being,  involves an abstract exploration of the concept of  'space' and the differential geometry used in mathematical physics. Hegel, who taught differential calculus both at high school and  university, dedicates a long section of the Logic to the infinitesimal calculus, the notions of which illuminate many other passages of the Logic. We start with what can be considered the archetype of the concept of space, the adjunct quadruple  associated to a the category of presheaves over a category $C$ with a terminal object. We call this space proto-space.

\[
\begin{tikzcd}
	PrShv(C) \arrow[r, "\Gamma"]  \arrow[r, shift left = 4ex, "\Pi"] & Set \arrow[l, "CoDisc", shift left = 2ex] \arrow[l, "Disc"', shift right = 2ex]
\end{tikzcd}
\]
Here $\Gamma(P) = P(1)$. It is very instructive to work out these adjoints $\Pi\vdash Disc \vdash\Gamma\vdash CoDisc$ explicitly. We have that $Disc (X)(A) = X$ for all objects $A$ in $C$ and $Disc(X)(f)$ is $id_X$ for any $ f : A \rightarrow B$ in $C$. We have that $\Pi (\mathcal{A}) := \bigcup_{U \in Obj C} \mathcal{A}(U) / \sim$ where $\sim$ is the equivalence relation which identifies $s$ and $s'$ over $U$ and $U'$ respectively if there are $f: V \rightarrow U$, $g : V \rightarrow U'$ with $s_V = s'_V$. If $\Omega$ is the subobject classifier for presheaves on a topological space $X$ then $\Pi(\Omega)$ gives the connected components of $X$ in the usual sense. We have finally that $CoDisc(X)(U) = \mathcal{P}X$. If we have a set map $f : \mathcal{A}(1) \rightarrow X$ then we
can define a morphisms of presheaves $f^\flat : \mathcal{A} \rightarrow Codisc(X)$ given by $f^\flat (U)(s)= \{ x \in X : \exists w \in \mathcal{A}(1), w_U = s \}$ for $s \in \mathcal{A}(U)$.
This quadruple is closely related to the above diagram for terminal and initial objects. The functors in the quadruple arise by taking the left and right Kan extensions along $\mathbb{T}$ and $1$. Thus we have the development of the concept of proto-space from that of being.
The terminal object is the limit of the empty diagram. (Co)limits and adjunctions are all special cases of Kan extensions. Thus starting from the 2-category of all categories we can derive the concept of proto-space employing only Kan extensions. Kan extensions express dialectical reason's process of passing to the other while preserving an essential mediating connection. Proto-space then assumes various form through various localisations giving rise to the sheaf toposes.

In general the adjoint quadruple does not carry over to a topos obtained by localisation. The condition of being cohesive is what guarantees this. One example of a cohesive topos are sheaves on a cohesive site. Thus we are lead spontaneously to the concept of cohesive topos as the right categorical notion of space.

But it should be stated that $Disc$ represents to movement towards discrete quantity or repulsion of units and $CoDisc$ represents the movement towards continuous quantity and coalescence of units.

The Yoneda embedding $C \rightarrow PrShv(C)$ expresses that each category $C$ unfolds into a proto-space. This unfolding of categories in the category of Being proceeds to something extrinsic, as a passage to the other. Dialectical reasoning asks for how something is constructed, for repressed history. A paradigm is that we can have diagrams into the category without the corresponding (co)limit existing.
But to think of a diagram and the concept of limit is already to posit the limit as something other and lacking in the category but nevertheless proceeding from it. Thus the complete and cocomplete category of presheaves given by the Yoneda embedding is a genuine Hegelian progression.

The key to understanding higher category theory is passing the above considerations into the correct generality of enriched category theory. A cosmos is what categories are enriched in. The above is the special case for the cosmos Set. Thus we should think of enriched presheaves, the enriched functor category between an enriched category and the cosmos itself seen as an enriched category.

Homotopy is the passage of quantity into quality. It is a changing of shape and size which preserves and thus defines a certain quality. Model categories are simply categorical abstractions where all constructions in classical homotopy theory can be carried out. The quality associated to a variation in quantity is expressed as the localisation yielding the homotopy category.

A simplicial set is an abstraction of a topological space, it is a categorical abstraction of a geometric form (i.e. a polytope). But it is a geometric form which contains within itself the process of its own genesis or assemblage, analogous to G-code (cf. the geometric realisation functor associating a topological space to each simplicial set). Category theory enriched over the cosmos of simplicial sets is currently seen to be the correct choice for the doing homotopy theory and differential geometry at the highest level of abstraction.

We must investigate the deeper significance of the simplex category $\Delta$ and its associated augmented simplex category $\Delta_a$ used to define simplicial sets. $\Delta_a$ has the natural structure of a strict monoidal category and $[0,1]$ has a natural monoid structure. This situation is universal in that monoidal categories $B$ with monoid objects $M$ are classified by functors $\Delta_a \rightarrow B$ sending $[0,1]$ to $M$. Similary we can obtain a classification of monads in a 2-category (curiously enough equivalent lax functors from the terminal object category). Perhaps $\Delta_a$ can be viewed as a higher qualitative categorical determination of the natural numbers expressing the exteriorisation and unfolding of the unit, one, monas (Hegel includes a short digression on Pythagoreanism in the section on Quantity). 

$Set$ is an exterior, abstract, discrete concept (Hegel's concept of number seems to be very set theoretic), but the category of simplicial sets $sSet$ represents a greater cohesion between parts and qualitative structure and determination. The morphisms between objects in a simplicially enriched category instead of being a mere set become a space.

To understand higher category theory it is important to master (among other things) the basics of the following subjects in category theory:

\begin{enumerate}

\item  Monoidal categories

\item  Enriched categories (including ends and coends and enriched Yoneda lemma)

\item  Model categories (including combinatorial model categories and homotopy (co)limits)  

\end{enumerate}

For 1 and 2 a good reference is Part I of Birgit Richter's From Categories to Homotopy Theory. For 3 it is the book chapter Homotopy theories and model categories by W. G. Dwyer and J. Spalinski.

One must feel very at home with simplicial sets (and their connection to homotopy types). Instead of ordinary (set-enriched) model categories we work with model categories enriched over simplicial sets sSet (seen itself as a category enriched over simplicial sets and when endowed with its classical (Quillen) model structure its fibrant-cofibrant objects are called $\infty$-groupoids).  Instead of presheaves to sets we work with simplicial set enriched functors from the opposite of a simplicial set enriched category to simplicial sets. A key result is Dugger's theorem, the simplicial presheaf analogue for combinatorial model categories of the characterization of sheaf toposes as left exact full subcategory localizations of presheaves - itself a generalization of the presentation of an object by generators and relations.  That it would seem that Topos Theory could be done (homotopically) over sSet rather than Set.  

It seems we can give a more geometric interpretation of the nerve of category $C$ (a canonical way of extracting an simplicial set ) given usually  in terms of composable sequences of arrows. For take three composable arrows $f,g,h$. Think of $f$ and $g$ as being in the plane but $h$ directed perpendicularly into space. Then we get 3-simplex in $N(C)_3$ with faces $(f,g, g\circ f)$, $(g,h, h\circ g)$ and $(g\circ f, h, h\circ (g \circ f))$ where we view a pair of composable arrows together with their composition as a triangle, i.e. a 2-simplex. Our intuition is that the nerve of a category keeps track of all commutative diagrams and each such diagram is a geometric object.

Of course there are competing definitions of (models of) $\infty$-groupoids besides the sSet-based one (Kan complexes) which itself is only one possible choice for shapes (which include cubical and cellular sets). All this suggests philosophically that our different concepts and models of what a 'space' is are special embodiments of a single 'pure' concept which is yet to be determined. Note that homotopy type theory connects $\infty$-groupoids to  (identity) types/ propositions or spaces of proofs/functions/computations.

\section*{Appendix 2: Aspects of Topos Theory and Sheaf Theory}
Is there a connection between the double-negation topology, the dense topology and the sheafication functor $P^{++}$ ?   The concept of sheaf is related to 'locality' and to the  'flowing' or 'variation' of one object relative to another. The locality aspect is seen  in Mitchel-Bénabou semantics, in the interpretation of $\vee$ and $\exists$ (reflected in their peculiar natural deduction rules) and the 'flowing' aspect can be appreciated as follows. When we define Dedekind reals in the internal language of a topos of sheaves over a topological space it turns out that these reals correspond to the sheaf of continuous functions over this topological space. So the sheaf condition has a subtle connection to continuity.

Consider the quotient category construction: we wish to construct a new category out of an old one in which a certain class $\Sigma$ of morphisms become isomorphisms. We are focusing on certain aspects of the objects ('equating' objects that share this aspect) and forgetting others. The category is seen according to as it is in a certain aspect. This is analogous to projection, for instance given, an idempotent linear operator. Now the following is  remarkable. From quotient categories we  are lead to the concept of a calculus of fractions and for this to work to the concept of reflection . But problems arise because the quotient construction can create more isomorphisms than is desired. Thus we require $r$ to be left exact (preserving finite limits) and we arrive at localisations. Under certain conditions this amounts to a closure operator and factorisation system.  Under certain conditions reflective subcategories correspond to idempotent monads (via their Eilenberg-Moore category). In type theory an idempotent monad defines a modality. Idempotent monads are the analogues of projection operators which focus on certain aspects of a being while negating others. We can order all (co) reflective subcategories of a category and subcategories which are both reflective and coreflective, the inclusion having both left and right adjoint. 

Four functors which form three adjunctions is the basis of Schreiber's theory of adjoint triples and the aufhebung relationship which is a pre-order having as bottom element the triple based on initial and terminal objects considered in a previous post.  In this pre-order there is a series of functors from solid to elastic to cohesive topoi and finally to set (which includes the global sections functor). 

Compactness is a generalisation of the concept of a finite set analogous to being finitely generated in algebra. Finitude or some control over 'size' is fundamental to define presentable and accessible categories, small objects, compact objects, generator, etc.). One idea is that any object can be given by a colimit of a diagram of generators. Many important categories are not small, their homsets are classes. So we need a way to control this largeness by having it determined in a some way by smallness, just as a the generators do for an algebraic structure. 

Consider the free completion or cocompletion $\mathcal{C} \rightarrow PShv(\mathcal{C})$. Absence of a (co)limit represents the 'bad infinity' which is overcome in the canonical completion (like the algebraic closure of a field). The presheaf category admits a further 'sheafication' functor $PShv(\mathcal{C}) \rightarrow Shv(\mathcal{C})$. This functor is a kind of localisation which dissolves global information (cf. quantum entanglement) and reduces it to information present in an  infinitesimal neighbourhood. Consider the case of the presheaf (which is not a sheaf) of constant functions on a topological space and observe what the sheafication functor does. If this presheaf were a sheaf then if we had two disconnected open subsets and if there were a section on the union of these subsets then the value on one would completely determine the value on another. For sheaves, the global is built from local information in a coherent way.

In Aristotle's \emph{Categories}  a division is made between discrete and continuous quantity.
If we visualize numbers - discrete quantity - in a spatial way, then spatiality is related to divisibility, to the formation of quotients involving equivalence relations.
If we take $0$ and $1$ then there is a space between them and we would like to name the "points" in this space and are thus
forced to "invent" $1/2$ with the additional caution that $1/2 \cong 2/4 \cong 4/8...$. This  abstract 
construction is often a stumbling block for children learning mathematics. 
Modern algebraic geometry springs from Descartes' theory of studying geometrical figures by means of algebraic equations (or polynomials).
We can study the behavior of varieties in a small vicinity of a  point - the topos of the point - analogously by
quotient constructions such as localisation and completion. 
There is also a similar concept of localization of a category. A full subcategory $i: D \rightarrow C$ is a localization if
there is a left adjoint $r$ to $i$ preserving limits. If furthermore $r$ is left exact then giving a localization is the same as
giving a \emph{calculus of fractions}, a construction analogous to rationals discussed above  where the morphisms of the category
are replaced by equivalence classes of certain pairs of morphisms $(f,g)$. 
We remarked how subobjects $Sub(A)$ of an object $A$ can be seen as a generalized space. It turns out
that giving a localization on a category is the same as giving a \emph{universal closure operator} on $Sub(A)$ for every object $A$.
This closure operator is a generalisation of the closure operator of ordinary topology. For a detailed discussion of the topics
of this section see \cite{Pedicchio}[Ch.7]. The inverse image sheaf functor we considered has a powerful expression
in the context of derived categories of Abelian sheaves and these are also given by a localization of
the homotopy category of complexes \cite{Kashiwara}[ch.1].  In fact modern homotopy theory is centered around the interplay between different kinds of simplicial structures and model categories\footnote{The notions of fibration and co-fibration can be seen as a generalization of what it means for a space to be respectively "inside" and "over" another space.}.  These "model" a homotopy theory which can be defined as a category with a class of morphisms -  weak equivalences - we wish to take as isomorphisms. We form the homotopy category by localising such categories at weak equivalences. Certain classes of models form themselves a category and have in turn a natural notion of weak equialence. Thus they can be seen as a "homotopy theory of homotopy theories"  \cite{Bergner} [p.82], surely the most complex evolution of the ancient concept of topos to date. We will return to homotopy theory in our final section about mereology.

The category $Set^{C^{op}}$ of presheaves over a category $C$, which is simply the category of functors from $C^{op}$
to $Set$ which are not required to satisfy any "gluing conditions",  can be seen as a sort of \emph{materia prima}.  Given a category $C$, endowing $C$ with a Grothendieck topology is the same as defining a universal closure operator on 
$Set^{C^{op}}$ and the associated localization is precisely the Grothendieck topos, the category of sheaves on $C$.

\end{document}